\documentclass[sn-mathphys,Numbered]{sn-jnl}

\usepackage{graphicx}%
\usepackage{multirow}%
\usepackage{amsmath,amssymb,amsfonts}%
\usepackage{amsthm}%
\usepackage{mathrsfs}%
\usepackage[title]{appendix}%
\usepackage{xcolor}%
\usepackage{textcomp}%
\usepackage{manyfoot}%
\usepackage{booktabs}%
\usepackage{algorithm}%
\usepackage{algorithmicx}%
\usepackage{afterpage}
\usepackage{algpseudocode}%
\usepackage{listings}%
\usepackage{tabularx}
\usepackage{bookmark}
\usepackage{longtable}
\newcolumntype{L}{>{\raggedright\arraybackslash}X}

\usepackage{anyfontsize}

\usepackage{xcolor}
\usepackage{pifont}
\usepackage{pbox}
\usepackage{stackengine}
\usepackage{scalerel}
\usepackage{xcolor,amssymb}
\newcommand{\cmark}{\textcolor{green}{\ding{51}}}
\newcommand{\xmark}{\textcolor{red}{\ding{54}}}
\newcommand\dangersign[1][1.8ex]{%
  \scaleto{\stackengine{0.3pt}{\scalebox{1.6}[1.4]{%
  \color{orange}$\blacktriangle$}}{\small\bfseries \textcolor{black}{!}}{O}{c}{F}{F}{L}}{#1}%
}

\newcommand{\appendixentry}[4]{
\vspace{0.1cm}
\noindent \textbf{Article #1}: #2
\\ \textit{Title}: #3
\\ \textit{Citations}: #4
}
\newcommand{\appendixdata}[8]{
\\ \textit{Fluid-related PDE(s)}: #1
\\ \textit{Primary outcome(s)}: #2
\\ \textit{Baseline}: #4
\\ \textit{Rule 1}: #5
\\ \textit{Rule 2}: #6
\\ \textit{Fair comparison}: #8
}

\newcommand{\fnoli}{1}
\newcommand{\fnolicitations}{941}
\newcommand{\citefnoli}{\citet{li2020fourier}}

\newcommand{\deeponetlulu}{2}
\newcommand{\deeponetlulucitations}{911}
\newcommand{\citedeeponetlulu}{\citet{lu2021learning}}

\newcommand{\tompson}{3}
\newcommand{\tompsoncitations}{558}
\newcommand{\citetompson}{\citet{tompson2017accelerating}}

\newcommand{\mlacceleratedcfd}{4}
\newcommand{\mlacceleratedcfdcitations}{429}
\newcommand{\citemlacceleratedcfd}{\citet{kochkov2021machine}}

\newcommand{\meshbasedpfaff}{5}
\newcommand{\meshbasedpfaffcitations}{390}
\newcommand{\citemeshbasedpfaff}{\citet{pfaff2020learning}}

\newcommand{\barsinai}{6}
\newcommand{\barsinaicitations}{382}
\newcommand{\citebarsinai}{\citet{bar2019learning}}

\newcommand{\deepfluids}{7}
\newcommand{\deepfluidscitations}{382}
\newcommand{\citedeepfluids}{\citet{kim2019deep}}

\newcommand{\deeponetwang}{8}
\newcommand{\deeponetwangcitations}{230}
\newcommand{\citedeeponetwang}{\citet{wang2021learning}}

\newcommand{\solverinthe}{9}
\newcommand{\solverinthecitations}{143}
\newcommand{\citesolverinthe}{\citet{um2020solver}}

\newcommand{\deepmmnetcvt}{10}
\newcommand{\deepmmnetcvtcitations}{128}
\newcommand{\citedeepmmnetcvt}{\citet{cai2021deepm}}

\newcommand{\belbuteperez}{11}
\newcommand{\belbuteperezcitations}{124}
\newcommand{\citebelbuteperez}{\citet{belbute2020combining}}

\newcommand{\pinoli}{12}
\newcommand{\pinolicitations}{124}
\newcommand{\citepinoli}{\citet{li2021physics}}

\newcommand{\projectionyang}{13}
\newcommand{\projectionyangcitations}{122}
\newcommand{\citeprojectionyang}{\citet{yang2016data}}

\newcommand{\neuralconverge}{14}
\newcommand{\neuralconvergecitations}{101}
\newcommand{\citeneuralconverge}{\citet{hsieh2019learning}}

\newcommand{\messagepassing}{15}
\newcommand{\messagepassingcitations}{87}
\newcommand{\citemessagepassing}{\citet{brandstetter2022message}}

\newcommand{\deepmmnet}{16}
\newcommand{\deepmmnetcitations}{84}
\newcommand{\citedeepmmnet}{\citet{mao2021deepm}}

\newcommand{\frameworkmishra}{17}
\newcommand{\frameworkmishracitations}{80}
\newcommand{\citeframeworkmishra}{\citet{mishra2018machine}}

\newcommand{\optimizemultigrid}{18}
\newcommand{\optimizemultigridcitations}{79}
\newcommand{\citeoptimizemultigrid}{\citet{greenfeld2019learning}}

\newcommand{\donga}{19}
\newcommand{\dongacitations}{55}
\newcommand{\citedonga}{\citet{dong2021local}}

\newcommand{\rayb}{20}
\newcommand{\raybcitations}{52}
\newcommand{\citerayb}{\citet{ray2019detecting}}

\newcommand{\novelcnn}{21}
\newcommand{\novelcnncitations}{47}
\newcommand{\citenovelcnn}{\citet{xiao2018novel}}

\newcommand{\wandel}{22}
\newcommand{\wandelcitations}{46}
\newcommand{\citewandel}{\citet{wandel2020learning}}

\newcommand{\shan}{23}
\newcommand{\shancitations}{46}
\newcommand{\citeshan}{\citet{shan2020study}}

\newcommand{\algebraic}{24}
\newcommand{\algebraiccitations}{46}
\newcommand{\citealgebraic}{\citet{luz2020learning}}

\newcommand{\zhuang}{25}
\newcommand{\zhuangcitations}{40}
\newcommand{\citezhuang}{\citet{zhuang2021learned}}

\newcommand{\pathak}{26}
\newcommand{\pathakcitations}{36}
\newcommand{\citepathak}{\citet{pathak2020using}}

\newcommand{\leoni}{27}
\newcommand{\leonicitations}{34}
\newcommand{\citeleoni}{\citet{di2021deeponet}}

\newcommand{\fnodeformation}{28}
\newcommand{\fnodeformationcitations}{33}
\newcommand{\citefnodeformation}{\citet{li2022fourier}}

\newcommand{\stevensa}{29}
\newcommand{\stevensacitations}{27}
\newcommand{\citestevensa}{\citet{stevens2020enhancement}}

\newcommand{\illarramendi}{30}
\newcommand{\illarramendicitations}{27}
\newcommand{\citeillarramendi}{\citet{ajuria2020towards}}

\newcommand{\stachenfeld}{31}
\newcommand{\stachenfeldcitations}{27}
\newcommand{\citestachenfeld}{\citet{stachenfeld2021learned}}

\newcommand{\han}{32}
\newcommand{\hancitations}{27}
\newcommand{\citehan}{\citet{han2022predicting}}

\newcommand{\stevensb}{33}
\newcommand{\stevensbcitations}{26}
\newcommand{\citestevensb}{\citet{stevens2020finitenet}}

\newcommand{\ozbay}{34}
\newcommand{\ozbaycitations}{25}
\newcommand{\citeozbay}{\citet{ozbay2021poisson}}

\newcommand{\zili}{35}
\newcommand{\zilicitations}{24}
\newcommand{\citezili}{\citet{li2022graph}}

\newcommand{\peng}{36}
\newcommand{\pengcitations}{23}
\newcommand{\citepeng}{\citet{peng2022attention}}

\newcommand{\chen}{37}
\newcommand{\chencitations}{21}
\newcommand{\citechen}{\citet{chen2021numerical}}

\newcommand{\alguacilb}{38}
\newcommand{\alguacilbcitations}{21}
\newcommand{\citealguacilb}{\citet{alguacil2021predicting}}

\newcommand{\wandelb}{39}
\newcommand{\wandelbcitations}{21}
\newcommand{\citewandelb}{\citet{wandel2021teaching}}

\newcommand{\blist}{40}
\newcommand{\blistcitations}{16}
\newcommand{\citeblist}{\citet{list2022learned}}

\newcommand{\cheng}{41}
\newcommand{\chengcitations}{15}
\newcommand{\citecheng}{\citet{cheng2021using}}

\newcommand{\wen}{42}
\newcommand{\wencitations}{13}
\newcommand{\citewen}{\citet{wen2020edge}}

\newcommand{\delaraa}{43}
\newcommand{\delaraacitations}{12}
\newcommand{\citedelaraa}{\citet{de2022accelerating}}

\newcommand{\zhao}{44}
\newcommand{\zhaocitations}{10}
\newcommand{\citezhao}{\citet{zhao2022learning}}

\newcommand{\assessments}{45}
\newcommand{\assessmentscitations}{8}
\newcommand{\citeassessments}{\citet{illarramendi2022performance}}

\newcommand{\holloway}{46}
\newcommand{\hollowaycitations}{7}
\newcommand{\citeholloway}{\citet{holloway2021acceleration}}

\newcommand{\azulay}{47}
\newcommand{\azulaycitations}{7}
\newcommand{\citeazulay}{\citet{azulay2022multigrid}}

\newcommand{\wulatent}{48}
\newcommand{\wulatentcitations}{7}
\newcommand{\citewulatent}{\citet{wu2022learning}}

\newcommand{\liu}{49}
\newcommand{\liucitations}{6}
\newcommand{\citeliu}{\citet{liu2022predicting}}

\newcommand{\zhang}{50}
\newcommand{\zhangcitations}{5}
\newcommand{\citezhang}{\citet{zhang2022hybrid}}

\newcommand{\duarte}{51}
\newcommand{\duartecitations}{4}
\newcommand{\citeduarte}{\citet{duarte2022black}}

\newcommand{\alguacil}{52}
\newcommand{\alguacilcitations}{4}
\newcommand{\citealguacil}{\citet{alguacil2022deep}}

\newcommand{\bezginb}{53}
\newcommand{\bezginbcitations}{4}
\newcommand{\citebezginb}{\citet{bezgin2022weno3}}

\newcommand{\shang}{54}
\newcommand{\shangcitations}{4}
\newcommand{\citeshang}{\citet{shang2022deep}}

\newcommand{\kube}{55}
\newcommand{\kubecitations}{3}
\newcommand{\citekube}{\citet{kube2021machine}}

\newcommand{\shi}{56}
\newcommand{\shicitations}{3}
\newcommand{\citeshi}{\citet{shi2022lordnet}}

\newcommand{\ranadea}{57}
\newcommand{\ranadeacitations}{3}
\newcommand{\citeranadea}{\citet{ranade2021latent}}

\newcommand{\chenb}{58}
\newcommand{\chenbcitations}{3}
\newcommand{\citechenb}{\citet{chen2022machine}}

\newcommand{\ranadeb}{59}
\newcommand{\ranadebcitations}{3}
\newcommand{\citeranadeb}{\citet{ranade2021composable}}

\newcommand{\pengb}{60}
\newcommand{\pengbcitations}{3}
\newcommand{\citepengb}{\citet{peng2023linear}}

\newcommand{\delarab}{61}
\newcommand{\delarabcitations}{2}
\newcommand{\citedelarab}{\citet{de2023accelerating}}

\newcommand{\ranadec}{62}
\newcommand{\ranadeccitations}{2}
\newcommand{\citeranadec}{\citet{ranade2022composable}}

\newcommand{\fang}{63}
\newcommand{\fangcitations}{2}
\newcommand{\citefang}{\citet{fang2022immersed}}

\newcommand{\shukla}{64}
\newcommand{\shuklacitations}{2}
\newcommand{\citeshukla}{\citet{shukla2023deep}}

\newcommand{\zhangb}{65}
\newcommand{\zhangbcitations}{2}
\newcommand{\citezhangb}{\citet{zhang2022learning}}

\newcommand{\bezgin}{66}
\newcommand{\bezgincitations}{2}
\newcommand{\citebezgin}{\citet{bezgin2021fully}}

\newcommand{\yang}{67}
\newcommand{\yangcitations}{2}
\newcommand{\citeyang}{\citet{yang2023rapid}}

\newcommand{\tang}{68}
\newcommand{\tangcitations}{2}
\newcommand{\citetang}{\citet{tang2022neural}}

\newcommand{\nastorg}{69}
\newcommand{\nastorgcitations}{1}
\newcommand{\citenastorg}{\citet{nastorg2022ds}}

\newcommand{\gopakumar}{70}
\newcommand{\gopakumarcitations}{1}
\newcommand{\citegopakumar}{\citet{gopakumar2023fourier}}

\newcommand{\shit}{71}
\newcommand{\shitcitations}{0}
\newcommand{\citeshit}{\citet{shit2021semi}}

\newcommand{\suforecast}{72}
\newcommand{\suforecastcitations}{0}
\newcommand{\citesuforecast}{\citet{su2021forecasting}}

\newcommand{\jeon}{73}
\newcommand{\jeoncitations}{0}
\newcommand{\citejeon}{\citet{jeon2022physics}}

\newcommand{\dai}{74}
\newcommand{\daicitations}{0}
\newcommand{\citedai}{\citet{dai2022fournetflows}}

\newcommand{\sun}{75}
\newcommand{\suncitations}{0}
\newcommand{\citesun}{\citet{sun2022local}}

\newcommand{\shao}{76}
\newcommand{\shaocitations}{0}
\newcommand{\citeshao}{\citet{shao2022poisson}}

\newcommand{\discacciati}{77}
\newcommand{\discacciaticitations}{54}
\newcommand{\citediscacciati}{\citet{discacciati2020controlling}}

\newcommand{\magiera}{78}
\newcommand{\magieracitations}{44}
\newcommand{\citemagiera}{\citet{magiera2020constraint}}

\newcommand{\bezginc}{79}
\newcommand{\bezginccitations}{16}
\newcommand{\citebezginc}{\citet{bezgin2021data}}

\newcommand{\dongb}{80}
\newcommand{\dongbcitations}{13}
\newcommand{\citedongb}{\citet{dong2022computing}}

\newcommand{\dresdner}{81}
\newcommand{\dresdnercitations}{8}
\newcommand{\citedresdner}{\citet{dresdner2022learning}}

\newcommand{\toshev}{82}
\newcommand{\toshevcitations}{0}
\newcommand{\citetoshev}{\citet{toshev20233}}

\begin{document}

\title[Weak baselines and reporting biases lead to overoptimism in machine learning for fluid-related partial differential equations]{Weak baselines and reporting biases lead to overoptimism in machine learning for fluid-related partial differential equations}

\author*[1,2]{\fnm{Nick} \sur{McGreivy}}\email{mcgreivy@princeton.edu}

\author[2]{\fnm{Ammar} \sur{Hakim}}\email{ahakim@pppl.gov}

\affil*[1]{\orgdiv{Department of Astrophysical Sciences}, \orgname{Princeton University}, \orgaddress{\city{Princeton}, \state{New Jersey}, \country{USA}}}
\affil[2]{\orgname{Princeton Plasma Physics Laboratory}, \orgaddress{\street{100 Stellarator Rd}, \city{Princeton}, \state{New Jersey}, \country{USA}}}

\abstract{One of the most promising applications of machine learning (ML) in computational physics is to accelerate the solution of partial differential equations (PDEs). The key objective of ML-based PDE solvers is to output a sufficiently accurate solution faster than standard numerical methods, which are used as a baseline comparison.
We first perform a systematic review of the ML-for-PDE solving literature.
Of articles that use ML to solve a fluid-related PDE and claim to outperform a standard numerical method, we determine that 79\% (60/76) compare to a weak baseline.
Second, we find evidence that reporting biases, especially outcome reporting bias and publication bias, are widespread.
We conclude that ML-for-PDE solving research is overoptimistic:
weak baselines lead to overly positive results, while reporting biases lead to underreporting of negative results.
To a large extent, these issues appear to be caused by factors similar to those of past reproducibility crises: researcher degrees of freedom and a bias towards positive results.
We call for bottom-up cultural changes to minimize biased reporting as well as top-down structural reforms intended to reduce perverse incentives for doing so.}

\keywords{machine learning, partial differential equations, metascience, reproducibility crisis}

\maketitle

\section{Introduction}\label{sec:introduction}
Many fields of science have experienced reproducibility issues \cite{randall2018irreproducibility,ritchie2020science,munafo2017manifesto}.
In some fields, reproducibility issues are thought to impact the validity of a significant percentage of published research \cite{ioannidis2005most,open2015estimating,prinz2011believe,begley2012raise}.
These issues are often caused by pitfalls with data analysis and statistical techniques, as well as by a systemic bias towards publishing positive results \cite{randall2018irreproducibility,ritchie2020science,munafo2017manifesto,gelman2013garden}.
Because these issues can undermine the credibility and authority of an entire field, they are often referred to as a `reproducibility crisis' \cite{baker2016reproducibility}.

As interest in machine learning (ML) has grown, more and more scientific fields are exploring whether ML can be used to advance science \cite{karagiorgi2022machine,dara2022machine,mater2019deep,carleo2019machine,rajkomar2019machine,grimmer2021machine}. 
For some problems, ML has shown the potential to do so \cite{jumper2021highly}. However, there are increasing concerns about reproducibility issues in ML \cite{hutson2018artificial,gundersen2022sources,sculley2018winner,armstrong2009improvements} and in ML-based science \cite{kapoor2023leakage,kapoor2023reforms}. Compiling evidence from 22 articles across 17 fields analyzing reproducibility issues in 294 articles, \citet{kapoor2023leakage} argue that there is a `reproducibility crisis' in ML-based science. 
Other large-scale analyses have found frequent reproducibility issues across hundreds of articles in medical ML \cite{demasi2017meaningless,roberts2021common,wynants2020prediction}.
Common pitfalls include data leakage \cite{kapoor2023leakage,whalen2022navigating}, poor data quality \cite{kapoor2023leakage,roberts2021common,artrith2021best}, weak baselines \cite{demasi2017meaningless}, and insufficient external validation \cite{roberts2021common,wynants2020prediction}. In each case, pitfalls result in overoptimistic assessments about the performance of ML.

In recent years, there has been interest in using ML to advance research into partial differential equations (PDEs) \cite{thuerey2021physics,brunton2023machine,vinuesa2022enhancing,karniadakis2021physics,cuomo2022scientific,duraisamy2019turbulence}. Scientists and engineers study PDEs because they accurately model the behavior of many physical systems.
PDEs relate the output of a function to partial derivatives with respect to input variables (usually space and/or time).
ML-for-PDE research has mostly focused on either solving well-posed `forward' initial boundary value problems, ill-posed `inverse' problems that use data to infer equation parameters or missing data \cite{karniadakis2021physics}, or `reduced order models' that learn low-dimensional representations from data \cite{brunton2023machine}.

Computational scientists have spent decades developing numerical algorithms to approximate the solution of forward PDE problems \cite{durran2013numerical,leveque1992numerical} -- we call these algorithms `standard numerical methods' or `standard solvers' -- but there is great interest in using ML to do so more efficiently \cite{brunton2023machine,vinuesa2022enhancing}.
Standard numerical methods have a basic trade-off between accuracy and efficiency; computing a more accurate approximation takes more time.
In principle, ML could be used to learn new algorithms or surrogate models that reduce the time required to output an approximate solution compared to standard solvers.
Faster ML-based solvers could be useful for downstream applications such as optimization, inverse problems, or uncertainty quantification \cite{mishra2018machine}, and to improve or even replace the standard numerical methods used in simulation codes for research and commercial applications \cite{kochkov2021machine}.
Indeed, many articles claim to have used ML to accelerate the solution of PDEs. These articles compare to standard numerical methods which serve as baselines and typically use data to train neural networks as components within surrogate models.

While various successes have been reported for ML-based PDE solvers compared to standard numerical methods, there are in our view several ways in which this success is likely limited.
First, to be useful as a surrogate model, an ML-based solver must reduce the total computational cost in downstream applications.
This includes the cost of generating data and training models, both of which are unaccounted for when comparing speed to standard solvers \cite{kadapa2021machine}.
Speed is thus necessary, but not sufficient, to be useful for forward problems.
Second, ML-based solvers have important qualitative limitations which standard solvers do not share. In particular, there are serious concerns about accuracy when generalizing to new parameter spaces \cite{ross2023benchmarking}, numerical stability \cite{lippe2024pde}, and predicting chaotic systems \cite{vlachas2020backpropagation}.
Increased speed seems to be the one area where ML-based solvers might have an advantage over standard solvers.
Third, one of the most widely researched methods involving ML and PDEs -- the so-called `physics-informed neural network (PINN) \cite{raissi2019physics} -- is known to be orders of magnitude slower than standard numerical methods at forward problems (i.e., solving PDEs) \cite{grossmann2023can,de2023physics}. Furthermore, PINNs often fail to converge to a reasonable approximation \cite{chuang2022experience,chuang2023predictive}, even for simple toy problems \cite{wang2022and,krishnapriyan2021characterizing,basir2022critical}. For ill-posed and inverse problems PINNs have had some success  \cite{karniadakis2021physics}, though they appear not to outperform alternatives such as discrete grid-based methods \cite{karnakov2024solving}.

The results in this Analysis call into question whether ML has actually been as successful at solving PDEs from fluid mechanics and related fields as the scientific literature would suggest. We identify two issues, weak baselines and reporting biases, that lead to overoptimism and affect the interpretation reproducibility \cite{gundersen2021fundamental} of ML-for-PDE solving research. To determine the frequency of weak baselines, we conduct a systematic review of research that uses ML to accelerate the solution of fluid-related PDEs. We observe two common pitfalls with baseline comparisons and identify them in a majority (79\%) of published articles. 
We then use anecdotal and statistical evidence from the systematic review to argue that reporting biases are causing negative results to be underreported.
Due to these reproducibility issues, we caution that at present the scientific literature is not a reliable source for evaluating the success of ML at solving PDEs.

\section{Weak baselines \label{sec:weakbaselines}}

ML-based PDE solvers use neural networks, deep learning, and other techniques from ML to output approximate solutions to forward PDE problems.
To determine whether ML can improve the efficiency of PDE solvers, ML-based solvers must be compared to standard numerical methods which serve as a baseline.
For these comparisons to reach accurate conclusions, they must be fair.
A fair comparison should not overestimate or underestimate the performance of either method.
Underestimating the performance of a baseline is called using a weak baseline. In ML-for-PDE solving research, we've observed two common pitfalls that lead to weak baselines. As a result, we introduce two rules (rule 1 and rule 2) that must be followed for comparisons between ML-based solvers and standard numerical methods to be fair. 
We also introduce three additional recommendations for fair comparisons (see Methods) but do not require that they be satisfied to ensure a fair comparison.

The first pitfall is to compare the efficiency of the (highly accurate) standard solver used to generate the training data to that of the (less accurate) ML-based solver.
The problem is that the standard solver can trade off accuracy for runtime, so it too can be made faster but less accurate.
If the ML-based solver has less accuracy, comparing the runtime with a standard solver is only meaningful if the standard solver is modified to compare at equal accuracy.
Comparing the runtime between two methods with different accuracies could lead, for example, to the nonsensical conclusion that a method is orders of magnitude faster than itself.
Thus, rule 1 is to make comparisons at either equal accuracy or equal runtime.
To satisfy rule 1, either (a) reduce the resolution and/or the number of iterations of the standard solver until the two methods have either equal accuracy or equal runtime (or a proxy for runtime), or (b) demonstrate that reducing the resolution of the standard solver any further would give worse accuracy than the ML-based solver.
This must be done even when the two solutions look qualitatively similar, because a lower resolution baseline will often look qualitatively similar as well.

The second pitfall is to compare to a numerical method which is much less efficient than a state-of-the-art method for that problem.
State-of-the-art numerical methods can be orders of magnitude faster than less efficient numerical methods.
Choosing the right algorithm for a given PDE can be difficult and requires a combination of expertise, background knowledge, and effort.
Thus, rule 2 is to compare to an efficient numerical method.
The methods that satisfy rule 2 depend on the particular PDE being solved; we discuss the criteria we used to evaluate rule 2 later in this section.

How frequently does the ML-for-PDE solving literature compare to weak baselines? To answer this question, we perform a systematic review \cite{aromataris2014systematic}. This systematic review attempts to find every article that (i) uses ML to solve a fluid-related PDE and (ii) compares speed, or some proxy for speed, to a standard numerical method.
We exclude articles that use the PINN method \cite{raissi2019physics} because it is already well-known that standard numerical methods outperform PINNs for forward problems \cite{grossmann2023can}.
We define additional inclusion and exclusion criteria in Methods.
We restrict ourselves to PDEs that are related to fluid mechanics both because these PDEs have received the most attention but also because this is our area of expertise. 

We find 82 articles (see Supplementary Information) matching the inclusion criteria.
76 articles claim to outperform and 4 claim to underperform relative to a standard numerical method \cite{magiera2020constraint,bezgin2021data,dresdner2022learning,toshev20233}, while 2 claim to have similar or varied performance \cite{discacciati2020controlling,dong2022computing}.
For each article that claims to outperform a standard numerical method, we ask whether that article's `primary outcome' (defined in Methods) followed rules 1 and 2.
We evaluate each rule using a three-point scale. We give a check ({\cmark}) if the rule was satisfied, or if we are unsure. We give an `X' ({\xmark}) if the rule was not satisfied. 
We give a warning sign ({\dangersign}) if the rule was only partially satisfied (rule 1) or if we believe the rule was likely not satisfied but we don't have enough evidence to say for sure (rule 2).
For rule 2, we gave an ({\xmark}) if
\begin{itemize}
    \item we could replicate the article’s primary outcome and achieve significantly improved performance with a different baseline (six articles).
    \item the article reported performance relative to a weaker baseline in the abstract but a stronger baseline in the results section or the appendix (three articles).
    \item the article used a 2D code as a baseline for a 1D problem (one article).
    \item the computational implementation of a state-of-the-art numerical method was orders of magnitude slower than our implementation of that method (one article).
    \item the baseline used implicit timestepping when explicit timestepping would have been faster (one article).
    \item for elliptic PDEs, the baseline method was much less efficient than a state-of-the-art numerical method (eight articles).
    \item for advection-dominated PDEs, the baseline method was much less efficient than a state-of-the-art numerical method (five articles).
\end{itemize}
For elliptic and advection-dominated PDEs, which have been the focus of our replication efforts, we explain which numerical methods we consider state-of-the-art in Methods.
Eight articles received a ({\dangersign}) for rule 2, either because (a) we believe the numerical method used is much less efficient than a state-of-the-art method, but we haven't performed a direct comparison for that PDE so we can't say for sure, or (b) the general-purpose software package being used as a baseline has been shown to be computationally slow relative to other solvers for that PDE, but we haven't performed direct comparisons ourselves and so we can't say for sure. One article received a ({\dangersign}) for rule 1 for reducing the resolution of the highly accurate numerical method, but not by enough to compare at equal accuracy.

We list the 76 articles that claim to outperform a standard numerical method in Table \ref{tab:baselines}, ordered from highest to lowest number of citations.
60/76 (79\%) receive an ({\xmark}) for rule 1 and/or rule 2 and thus compare to a weak baseline. 2/76 (2.6\%) receive a ({\dangersign}), indicating that they may be comparing to a weak baseline and the performance claims should be treated with caution. 14/76 (18.4\%) receive a ({\cmark}), indicating that we believe they compare to a strong baseline.
See Supplementary Information for detailed explanations of each entry in Table \ref{tab:baselines}.
Articles which receive a ({\cmark}) tend to have quantitatively smaller relative improvements than articles which receive an ({\xmark}), suggesting that the more impressive the result, the more likely the article used a weak baseline. 

\setlength\LTcapwidth{1.1\linewidth}
\setlength{\tabcolsep}{3pt}
\begin{longtable}[p!]{@{\extracolsep\fill}llllccc}
\caption{Weak baselines in ML-for-PDE solving research.}\\
\endfirsthead
\toprule
\endhead
\bottomrule
\endfoot
\endlastfoot
\toprule
Article & Cited & PDE & Primary outcome & Rule 1 & Rule 2 & Fair? \\
\midrule

    \fnoli{ }\citefnoli& \fnolicitations & a,b &up to $10^3 \times$ faster& \xmark & \cmark & \xmark \\
    \deeponetlulu{ }\citedeeponetlulu& \deeponetlulucitations & c,d& substantially lower cost & \xmark & \cmark & \xmark \\
    \tompson{ }\citetompson& \tompsoncitations & e &$4.1\times$ faster& \cmark & \cmark & \cmark \\
    \mlacceleratedcfd{ }\citemlacceleratedcfd& \mlacceleratedcfdcitations & b & 8-10$\times$ coarser, 40-80$\times$ speedup & \cmark & \xmark & \xmark \\
    \meshbasedpfaff{ }\citemeshbasedpfaff& \meshbasedpfaffcitations & f,g & $10^1$-$10^2\times$ faster& \xmark & \cmark& \xmark \\
    \barsinai{ }\citebarsinai& \barsinaicitations & a,h,i &4-8$\times$ coarser& \cmark & \xmark & \xmark \\
    \deepfluids{ }\citedeepfluids& \deepfluidscitations & j &700$\times$ faster& \xmark & \cmark & \xmark \\
    \deeponetwang{ }\citedeeponetwang& \deeponetwangcitations & a &up to $10^3\times$ faster& \xmark & \cmark & \xmark \\
    \solverinthe{ }\citesolverinthe& \solverinthecitations &a,b,e,f&68$\times$ speedup& \xmark & \cmark & \xmark  \\
    \deepmmnetcvt{ }\citedeepmmnetcvt & \deepmmnetcvtcitations & k &$10^4\times$ speedup& \xmark & \cmark& \xmark \\
    \belbuteperez{ }\citebelbuteperez & \belbuteperezcitations & g & substantial speedup & \xmark & \cmark & \xmark \\
    \pinoli{ }\citepinoli& \pinolicitations & a,b & orders of magnitude speedup & \cmark & \xmark & \xmark \\
    \projectionyang{ }\citeprojectionyang& \projectionyangcitations & e & drastic speedup & \xmark & \xmark & \xmark \\
    \neuralconverge{ }\citeneuralconverge& \neuralconvergecitations & e &2-3$\times$ speedup& \cmark & \xmark & \xmark \\
    \messagepassing{ }\citemessagepassing& \messagepassingcitations & a,j,l &outperforms state-of-the-art & \cmark & \xmark & \xmark \\
    \deepmmnet{ }\citedeepmmnet&\deepmmnetcitations & m & over $10^5\times$ faster & \xmark & \xmark & \xmark \\
    \frameworkmishra{ }\citeframeworkmishra& \frameworkmishracitations & a,c,n &significant gain in efficiency& \cmark & \xmark & \xmark \\
    \optimizemultigrid{ }\citeoptimizemultigrid& \optimizemultigridcitations & o & improved convergence rates & \cmark& \cmark & \cmark \\
    \donga{ }\citedonga& \dongacitations & a,c,e,t & often exceeds performance & \cmark & \xmark & \xmark \\
    \rayb{ }\citerayb& \raybcitations & a,c,n &outperforms TVB limiter& \cmark & \cmark & \cmark \\
    \novelcnn{ }\citenovelcnn& \novelcnncitations & e & $10^2\times$ faster projection step & \xmark & \xmark& \xmark \\
    \wandel{ }\citewandel& \wandelcitations & f & $11$/$40\times$ faster on CPU/GPU & \xmark & \cmark & \xmark \\
    \shan{ }\citeshan& \shancitations & e & significant speedup & \xmark & \xmark & \xmark \\
    \algebraic{ }\citealgebraic& \algebraiccitations & o & improved convergence rates & \cmark & \cmark & \cmark \\
    \zhuang{ }\citezhuang& \zhuangcitations & c & 4$\times$ lower resolution & \cmark & \xmark & \xmark \\
    \pathak{ }\citepathak& \pathakcitations & b & lower resolution & \cmark & \cmark & \cmark \\
    \leoni{ }\citeleoni& \leonicitations & p & very small cost & \xmark & \cmark & \xmark \\
    \fnodeformation{ }\citefnodeformation& \fnodeformationcitations & b,g & $10^5\times$ faster & \xmark & \cmark & \xmark \\
    \stevensa{ }\citestevensa& \stevensacitations & a,c,n &outperforms WENO& \cmark & \cmark & \cmark \\
    \illarramendi{ }\citeillarramendi& \illarramendicitations & e & 3.2$\times$ faster & \cmark & \xmark & \xmark \\
    \stachenfeld{ }\citestachenfeld& \stachenfeldcitations & b,q & outperforms state-of-the-art & \cmark & \dangersign & \dangersign \\
    \han{ }\citehan& \hancitations & f,q & 100-800$\times$ speedup & \xmark & \dangersign & \xmark \\
    \stevensb{ }\citestevensb& \stevensbcitations & a,c,i & 2-3$\times$ lower error & \cmark & \cmark & \cmark \\
    \ozbay{ }\citeozbay& \ozbaycitations & e & improved preconditioner & \cmark & \cmark & \cmark \\
    \zili{ }\citezili& \zilicitations & j & 5-8$\times$ acceleration & \cmark & \xmark & \xmark \\
    \peng{ }\citepeng& \pengcitations & b & 8000$\times$ speedup & \xmark & \xmark & \xmark \\
    \chen{ }\citechen& \chencitations & b & 300-600$\times$ speedup & \xmark & \dangersign & \xmark \\
    \alguacilb{ }\citealguacilb& \alguacilbcitations & l & 15.5$\times$ speedup & \xmark & \dangersign & \xmark \\
    \wandelb{ }\citewandelb& \wandelbcitations & f &  considerably faster & \xmark &  \cmark & \xmark \\
    \blist{ }\citeblist& \blistcitations & b & 14.4$\times$ speedup & \cmark & \xmark & \xmark \\
    \cheng{ }\citecheng& \chengcitations & e & 2$\times$ faster & \cmark & \cmark & \cmark \\
    \wen{ }\citewen& \wencitations & n,r & fewer grid points & \cmark & \cmark & \cmark \\
    \delaraa{ }\citedelaraa& \delaraacitations & a & significant cost savings & \xmark & \cmark & \xmark \\
    \zhao{ }\citezhao& \zhaocitations & f,l & 8$\times$ or $35\times$ faster & \dangersign & \cmark & \dangersign \\
    \assessments{ }\citeassessments& \assessmentscitations & e & 10-25$\times$ faster & \cmark & \xmark & \xmark \\
    \holloway{ }\citeholloway& \hollowaycitations & s & $270\times$ speedup & \xmark& \cmark& \xmark \\
    \azulay{ }\citeazulay& \azulaycitations & t & favorable runtime on GPU & \cmark & \cmark & \cmark \\
    \wulatent{ }\citewulatent& \wulatentcitations & a,b & 840$\times$ speedup & \xmark & \cmark & \xmark \\
    \liu{ }\citeliu& \liucitations & a,b & 10-60$\times$ speedup & \xmark & \xmark & \xmark \\
    \zhang{ }\citezhang& \zhangcitations & e,t & up to $\mathcal{O}(10^2)$ more efficient & \cmark & \xmark & \xmark \\
    \duarte{ }\citeduarte& \duartecitations & u & $10^4\times$ faster & \xmark & \cmark & \xmark \\
    \alguacil{ }\citealguacil& \alguacilcitations & l & 141$\times$ acceleration & \xmark & \dangersign & \xmark \\
    \bezginb{ }\citebezginb& \bezginbcitations & c,n & similar or better performance & \cmark & \xmark & \xmark \\
    \shang{ }\citeshang& \shangcitations & e,l & much more accurate & \cmark & \cmark & \cmark \\
    \kube{ }\citekube& \kubecitations & v & 25\% fewer iterations & \xmark & \cmark & \xmark \\
    \shi{ }\citeshi& \shicitations & b,e & over $50\times$ faster & \xmark & \cmark & \xmark \\
    \ranadea{ }\citeranadea& \ranadeacitations & w & over 200$\times$ speedup & \xmark & \cmark & \xmark \\
    \chenb{ }\citechenb& \chenbcitations & e & fewer iterations & \cmark & \xmark & \xmark \\
    \ranadeb{ }\citeranadeb& \ranadebcitations & b,w,x & 40-100$\times$ faster & \xmark & \cmark & \xmark \\
    \pengb{ }\citepengb& \pengbcitations & b & 20$\times$ speedup & \xmark & \cmark & \xmark \\
    \delarab{ }\citedelarab& \delarabcitations & q & 4-5$\times$ faster & \cmark & \cmark & \cmark \\
    \ranadec{ }\citeranadec& \ranadeccitations & e,w,x & 40-50$\times$ faster & \xmark & \cmark & \xmark \\
    \fang{ }\citefang& \fangcitations & f & 38.5\% faster & \xmark & \cmark & \xmark \\
    \shukla{ }\citeshukla& \shuklacitations & g & 32,253$\times$ speedup & \xmark & \cmark & \xmark \\
    \zhangb{ }\citezhangb& \zhangbcitations & l & around $10\times$ faster & \xmark & \cmark & \xmark \\
    \bezgin{ }\citebezgin& \bezgincitations & q & outperforms Rusanov flux & \cmark & \xmark & \xmark \\
    \yang{ }\citeyang& \yangcitations & l & nearly $10^2\times$ faster & \xmark & \cmark & \xmark \\
    \tang{ }\citetang& \tangcitations & e & up to $12\times$ speedup & \cmark & \xmark & \xmark \\
    \nastorg{ }\citenastorg& \nastorgcitations & e & $10\times$ faster & \xmark & \cmark & \xmark \\
    \gopakumar{ }\citegopakumar& \gopakumarcitations & y & $10^6\times$ faster & \xmark & \dangersign & \xmark \\
    \shit{ }\citeshit& \shitcitations & d & 19.2\% faster & \cmark & \xmark & \xmark \\
    \suforecast{ }\citesuforecast& \suforecastcitations & q & over $10^3\times$ faster & \xmark & \cmark & \xmark \\
    \jeon{ }\citejeon& \jeoncitations & b & 1.8$\times$ acceleration & \xmark & \dangersign & \xmark \\
    \dai{ }\citedai& \daicitations & f & $10^4$-$10^5\times$ faster & \xmark & \dangersign & \xmark \\
    \sun{ }\citesun& \suncitations & e,t & better accuracy & \cmark & \cmark & \cmark \\
    \shao{ }\citeshao& \shaocitations & e & improves convergence & \cmark & \xmark & \xmark \\

\bottomrule
\caption*{\small{PDEs: (a) Burgers' (b) incompressible Navier-Stokes (INS) (c) advection (d) advection-diffusion (e) Poisson (f) INS wake dynamics (g) compressible Navier-Stokes airfoil wing (h) Korteweg-de Vries (i) Kuramoto-Sivashinsky (j) INS graphics/particle-based (k) electroconvection (l) wave (m) reacting Navier-Stokes (n) Euler (o) elliptic diffusion (p) parabolized stability equations (q) compressible Navier-Stokes (r) shallow water (s) Boltzmann collision operator (t) Helmholtz (u) black hole hydrodynamics (v) particle-in-cell (w) convective heat transfer (x) Laplace (y) magnetohydrodynamics}}
%\end{tabular*}
\label{tab:baselines}
\end{longtable}

\begin{table}[ht]
\caption{Reproducing results in ML-for-PDE solving research using stronger baselines.}\label{tab:reproducebaselines}%
\setlength{\tabcolsep}{2pt}
\begin{tabular}{@{}lllp{22mm}p{20mm}p{23mm}p{23mm}@{}}
\toprule
Article & Cited  & PDE & Weaker \newline baseline & Stronger \newline baseline & Old \newline outcome & New \newline outcome \\
\midrule
\fnoli{ }\citefnoli & \fnolicitations & b & PS $64\times64$ & DG2 $7\times7$ & $10^3\times$ faster & $7\times$ faster \\
\deeponetlulu{ }\citedeeponetlulu & \deeponetlulucitations & c  & FD $n_x$=100 & DG2 $n_x$=13 & $24\times$ faster & $10\times$ slower \\
\mlacceleratedcfd{ }\citemlacceleratedcfd & \mlacceleratedcfdcitations & b & FV & PS & 80$\times$ faster & slightly slower$^\dagger$ \\
\barsinai{ }\citebarsinai& \barsinaicitations & a & WENO & DG2/DG3 & 4-8$\times$ fewer DOF & 2-4$\times$ fewer DOF\\
\deeponetwang{ }\citedeeponetwang& \deeponetwangcitations & a & SP $n_x$=100 & FV $n_x$=100 & $10^3\times$ faster & 10$\times$ slower \\
\pinoli{ }\citepinoli& \pinolicitations & b & PS $64\times64$ & DG2 $3\times3$ & $10^3\times$ faster & $7\times$ faster \\
\neuralconverge{ }\citeneuralconverge& \neuralconvergecitations & e  & MG & LU & faster & $10^3\times$ slower \\
\messagepassing{ }\citemessagepassing& \messagepassingcitations & a,l & WENO, PS & WENO, FV & much faster & $10^3\times$ slower \\
\delaraa{ }\citedelaraa& \delaraacitations & a & DG28 $n_x$=1 & DG9 $n_x=1$ & 22-75$\times$ faster & 4-10$\times$ slower \\
\tang{ }\citetang& \tangcitations & e & CG \& MG & LU & $12\times$ faster & 35-500$\times$ slower \\
\botrule
\end{tabular}
\footnotetext{PDEs: (a) 1D Burgers' (b) 2D incompressible Navier-Stokes (c) 1D advection (e) 2D Poisson (l) 1D wave \newline Numerical methods: (PS) pseudo-spectral (FD) second-order finite-difference (DG2) Discontinuous galerkin, polynomial basis functions of order 2 (FV) finite volume (SP) spectral (MG) multigrid (LU) LU decomposition (CG) conjugate gradient \newline Abbreviations: (DOF) degrees of freedom ($n_x$) number of cells in the $x$ direction ($64\times64$) 64 cells in both the $x$ and $y$ directions }
\footnotetext[\dagger]{This result on GPU is consistent with \citet{dresdner2022learning} who reproduce the result on TPU.}
\end{table}

We reproduce results from ten articles using a stronger baseline; these articles are listed in Table \ref{tab:reproducebaselines}.
We primarily focus on reproducing results from highly cited articles solving 1D and 2D PDEs in regular geometries.
4/10 articles violate rule 1, while at least 6/10 violate rule 2 because they use an inefficient numerical method and/or an inefficient implementation.
In 9/10 cases the stronger baseline is at least two orders of magnitude faster than the slower baseline; the exception is article 6 \cite{bar2019learning}.
In 7/10 cases the stronger baseline outperforms the ML-based solver; the exceptions are articles 1 \cite{li2020fourier}, 6 \cite{bar2019learning}, and 12 \cite{li2021physics}. In Methods, we give additional details about each of the stronger baselines listed in Table \ref{tab:reproducebaselines}. We give additional details about the results of each reproducibility experiment in Supplementary Information.

\section{Reporting biases \label{sec:reportingbiases}}
Reporting biases is an umbrella term for a set of biases that can arise when the analysis, reporting, and/or interpretation of research findings are influenced by the nature and direction of results \cite{reportingbiascatalogue}.
Types of reporting biases include publication bias \cite{thornton2000publication}, spin bias \cite{boutron2018misrepresentation}, and different flavors of outcome reporting bias \cite{outcomereportingbiascatalogue} such as p-hacking \cite{head2015extent}, selective reporting \cite{saini2014selective}, outcome switching \cite{altman2017outcomeswitching}, and data-dredging \cite{erasmus2022dredging}.

Because reporting biases cause negative results to get suppressed \cite{de2018cumulative}, the percentage of positive results is believed to correlate with the frequency of reporting biases \cite{ritchie2020science,fanelli2010positive}.
To estimate the percentage of positive results in ML-for-PDE solving research, we analyzed a random sample of articles (see Methods). Of articles whose abstracts mention positive and/or negative experimental results, 94.8\% (220/232) mention only positive results, 5.2\% (12/232) mention both positive and negative results, and 0\% (0/232) mention only negative results.
This is an unusually high percentage of positive results compared to other fields of science \cite{fanelli2010positive} and motivates us to investigate whether reporting biases are causing negative results to be underreported in ML-for-PDE research.

During our systematic review, we found anecdotal and statistical evidence of publication bias and outcome reporting bias.
Of the 82 articles matching the inclusion criteria, 76 (93\%) claimed to outperform a standard numerical method baseline, while only 4 (5\%) claim to underperform relative to a baseline.
This suggests that many negative results are either not being published (publication bias) or being reported as positive due to outcome reporting bias.
Sure enough, a close reading of articles in the systematic review reveals evidence of outcome reporting bias, especially selective reporting and outcome switching: reporting the runtime of some PDEs but not others \cite{li2020fourier,wang2021learning,um2020solver,brandstetter2022message,dong2021local,shang2022deep}, only highlighting the results from the most successful PDE \cite{li2021physics,li2022fourier}, reporting performance in non-standard ways to seem more successful or to conceal a negative result \cite{bar2019learning,zhuang2021learned,kube2021machine,wang2021learning,stevens2020finitenet,dresdner2022learning,alguacil2021predicting,alguacil2022deep,shang2022deep}, or comparing to a stronger baseline in the main text but a weaker baseline in the abstract \cite{dong2021local,bezgin2022weno3,xiao2018novel}.
By cross-referencing with other articles, we also find evidence consistent with publication bias: some methods, which are successful on one PDE \cite{sanchez2020learning,wang2021learning,li2020fourier,bar2019learning}, either have worse performance when tested with different parameters or on different PDEs \cite{klimesch2022simulating,wang2022improved,gupta2022towards} or don't reproduce under nearly identical conditions \cite{mcgreivy2023invariant}.

We can more directly observe the collective effects of outcome reporting bias using a natural experiment in the ML-for-PDE solving literature. 
We gather two groups of articles, which we call sample A and sample B. 
Sample A includes all 76 articles in Table \ref{tab:baselines}.
Sample B is a random sample of 60 articles which use the PINN method to solve a fluid-related PDE.
There is one key difference between sample A and sample B: while every article in sample A claims to outperform a standard numerical method in speed or computational cost, the authors of every article in sample B believe that their ML-based solver underperforms relative to standard numerical methods
(see Methods for an explanation of why we can assume that they believe this).
If outcome reporting bias were not present, we would expect -- given that both samples try to solve fluid-related PDEs using ML and both report the accuracy of their proposed method relative to standard solvers -- that both samples would report the efficiency relative to standard solvers at similar rates.
Yet the percentage of articles that report this in the abstract is 80\% (61/76) in sample A and 8\% (5/60) in sample B. Only 12\% (7/60) of articles in sample B report the efficiency relative to standard solvers in the entire article.
In other words, when articles have a positive result they almost always highlight it, but when they have a negative result they rarely report it. 

\begin{figure}[ht]
\centering
\includegraphics[width=0.65\textwidth]{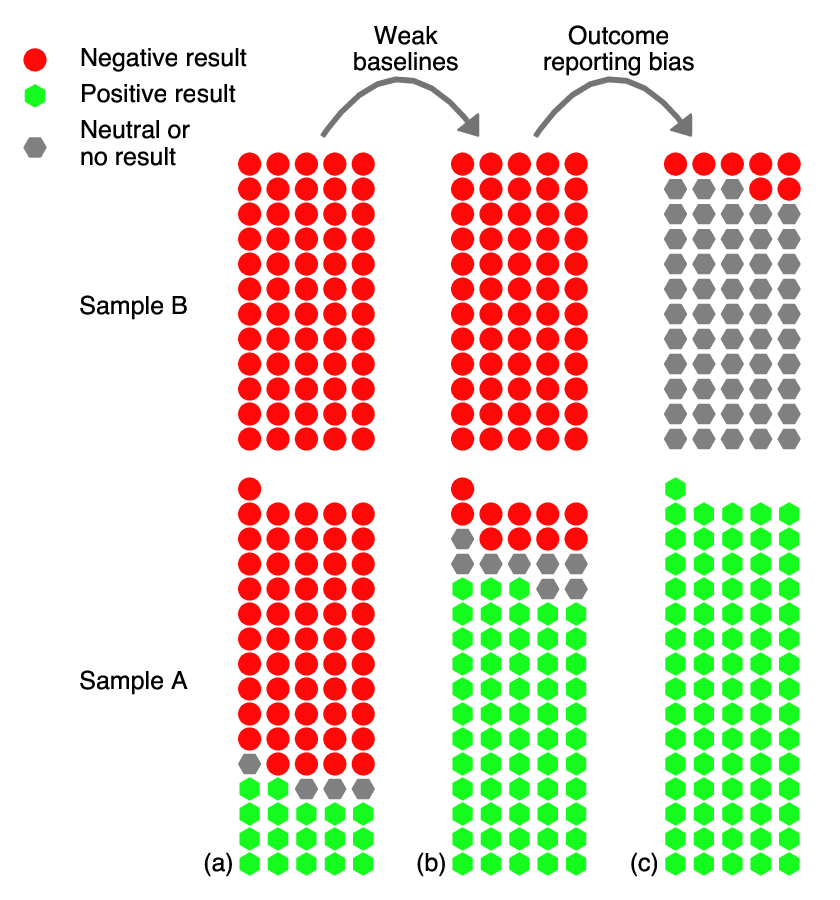}
\caption{The cumulative effects of weak baselines and reporting biases on samples A and B. Each circle or hexagon represents an article, while each color represents the result of comparing the relative speed and accuracy to a standard numerical method. In (a) we estimate what the results would be with strong baselines and without outcome reporting bias. (b) shows what the results would likely be without outcome reporting bias. (c) shows the results in the published literature.}
\label{fig:cumulativeeffects}
\end{figure}

Figure \ref{fig:cumulativeeffects} shows the cumulative effects of weak baselines and outcome reporting bias on samples A and B.
Weak baselines lead to overly positive results, while reporting biases lead to underreporting of negative results. The result is overoptimism about ML.

\section{Limitations
\label{sec:limitations}}
This Analysis has three main limitations.
First, the systematic review only considered forward problems and PDEs related to fluid mechanics.
While this is a large research area, much of the research on ML for scientific simulation is not analyzed (see exclusion criteria in Methods).
It is unclear whether these reproducibility issues also affect research using ML for inverse problems, solid mechanics PDEs, quantum mechanics PDEs, and high-dimensional PDEs.
Second, we do not have proof -- only evidence -- that reporting biases are causing negative results to be underreported.
It is possible that other factors could be contributing to the evidence we've presented, though we aren't sure what those factors might be.
Third, the bottom left column in Figure \ref{fig:cumulativeeffects} is an estimate.
There is also some uncertainty in the bottom center column, though less so.
The purpose of Figure \ref{fig:cumulativeeffects} is to illustrate the overall trend; readers should not conclude that we are making precise claims about what the results would be with stronger baselines and complete reporting.
%Fourth, discussion about the causes of these issues is necessarily speculative. Weak baselines and reporting biases are complex issues without monocausal explanations, and so there may be additional causes besides those we discuss below.

\section{Discussion \label{sec:conclusion}} 

To some extent, the issue of weak baselines (especially violations of rule 2) appears to be caused by three factors specific to this subfield: a lack of domain expertise in the ML community, insufficient benchmarking by the numerical analysis community, and little awareness about the difficulty and importance of choosing a strong baseline.
While benchmarking is a critical way of evaluating model performance in ML research, numerical analysis research tends to focus more on the theoretical properties of algorithms.
Furthermore, the relative performance of different numerical methods depends heavily on PDE parameters and implementation details.
It can be quite difficult, even for researchers with years of experience developing algorithms to solve PDEs, to predict the relative performance of different numerical methods on a given PDE.
Although a few researchers have begun to informally discuss issues with baseline comparisons in ML-for-PDE research, widespread awareness about the extent of the problem appears to be missing.

To reduce the frequency of baseline-related reproducibility issues in ML-for-PDE solving research, we make a few recommendations for best practices.
While failing to beat a baseline should not cause an article to be rejected, failing to follow best practices when evaluating models can and should be treated as grounds for rejection.
First, we recommend that all articles using ML to solve PDEs compare to two types of baselines: standard numerical methods and other ML-based solvers.
This allows readers to better evaluate model performance and reduces selective reporting.
A good example is \citet{stachenfeld2021learned}.
If other ML-based solvers cannot be implemented as baselines, articles should explain why.
Second, articles should follow rule 1 when comparing to standard numerical methods.
To satisfy rule 1, reduce the spatial resolution and/or the number of iterations of the standard solver until the two methods have equal accuracy or equal runtime.
Ideally, articles would also make plots of cost versus accuracy.
A good example is figure 1a of \citet{kochkov2021machine}.
Third, articles should discuss in a separate paragraph or subsection how the baselines were chosen and justify why the comparisons are unbiased. In particular, articles should explain why the standard numerical method being used as a baseline is highly efficient (or state-of-the-art) for that PDE.
Ideally, articles would compare multiple numerical methods for each PDE and use the most efficient method as a baseline.
A good example is Appendix C of \citet{cheng2021using}.
If authors are unsure which baseline is state-of-the-art for a given PDE, they should (a) talk to and/or collaborate with domain scientists or other experts, and (b) clearly acknowledge their uncertainty. A good example of (b) is the speed comparison in appendix D.1 of \citet{lippe2024pde}.
Fourth, besides rules 1 and 2 we make three additional recommendations for fair comparisons (see Methods).

To a large extent, however, weak baselines (especially violations of rule 1) and reporting biases appear be caused by factors similar to those that have led to reproducibility crises in other fields: researcher degrees of freedom, combined with a bias towards positive results.
In the process of writing an article about ML for PDEs, researchers make many choices.
Researchers choose not only PDEs, boundary conditions, hyperparameters, evaluation metrics, etc., but also which hypotheses to test, experiments to report, and results to emphasize.
Choices about experimental design, analysis, and reporting are called researcher degrees of freedom \cite{simmons2011false,wicherts2016degrees}.
Researcher freedom is valuable, but it becomes a problem when decisions about analysis and reporting are made or altered after results are known \cite{gelman2013garden}.
The conditional probability distribution of each decision given the experimental results tends to be biased in favor of positive outcomes.
The cumulative sum of these many biased decisions can significantly affect the reproducibility of published research \cite{de2018cumulative}.

We emphasize both good intentions and perverse incentives as explanations for the apparent bias towards positive results.
The culture of scientific ML is one in which well-intentioned researchers try to figure out ways that ML might be useful for science.
In the process of doing so, they tend to be less interested in reporting ways that ML isn't useful. 
Perverse incentives also contribute.
Because ML research rewards novel ideas and positive experimental results \cite{sculley2018winner}, all else being equal articles with weak baselines and/or reporting biases are more likely to get accepted to prestigious venues and more likely to be widely cited \cite{serra2021nonreplicable}.
Incentives against negative results are particularly strong in scientific ML, because career advancement (in academia) and lucrative jobs (in industry) depend on the presumption that ML will be a useful tool for scientific problems.
Negative results could cast doubt on that presumption, thereby undermining justification for one's research area.

Regardless of whether negative results are valuable in ML research \cite{borji2018negative}, the fact is that ML is now being used in science. In science, negative results matter. Without negative results, scientists cannot accurately determine whether and how ML is useful for advancing knowledge in their field. Unchecked, overoptimism can lead to misunderstanding of applicability, misallocation of research priorities, and slowdown in scientific progress.

Because the causes of biased reporting in ML-for-PDE research appear to be similar to those of past reproducibility crises, we recommend two types of reforms similar to those implemented by other fields: bottom-up cultural changes intended to minimize biased reporting, as well as top-down structural reforms intended to reduce perverse incentives for doing so. Most of these reforms will not benefit ML researchers – research projects will require more work and report more negative results -  but they will benefit science.

Cultural changes start at the level of the individual, research group, and/or department. ML researchers should have widespread awareness about and understanding of reproducibility issues. \citet{ritchie2020science} is a good place to start, while \citet{gundersen2022sources} discusses many issues unique to ML. Research groups should develop internal controls to ensure that reporting is complete and unbiased. Departments should teach about reproducibility pitfalls in ML classes. Individuals should commit to eliminating biased reporting, using strong baselines, discussing limitations honestly and transparently \cite{smith2022real}, and publishing negative results. To increase confidence in their conclusions, researchers can preregister their experiments or use registered reports.

We recommend two structural reforms.
First, ML journals and conferences could allow for registered reports \cite{gundersen2021case}.
These would be peer-reviewed before experiments are performed and evaluated based on whether the proposed analysis answers an interesting research question and is methodologically sound.
Accepted proposals would be guaranteed publication in the journal or a future conference, so long as the final paper conforms to the registered report.
Second, funding agencies could fund domain scientists to propose and setup challenge problems for the ML community to tackle.
A desirable challenge problem for scientific ML -- for example the CASP protein folding challenge \cite{jumper2021highly} -- would have three qualities.
First, the problem must be unsolved or extremely tedious to solve using standard methods.
Second, there must be a way of verifying whether or how well the problem was solved.
Third, the scientific community must agree that solving the problem would be a valuable contribution to science.
A challenge problem with these qualities would have clear evaluation metrics and either a standard baseline or no baseline at all, thereby eliminating the potential for weak baselines, reducing opportunities for outcome reporting bias, while also directing ML research away from toy problems and towards those of greatest scientific importance.

\section{Methods \label{sec:methods}}
We use the scientific literature to study issues with the scientific literature. This Analysis can thus be understood as an example of metascience, a field of research related to the use of scientific methods to study and improve science \cite{schooler2014metascience}.

\subsection{Systematic review}

A systematic review attempts to answer a predefined research question by collecting and analyzing evidence from all available research studies on the topic. Our research question is: how frequently does the ML-for-PDE solving literature compare to weak baselines?

\subsubsection{Inclusion criteria}

To restrict the scope of the systematic review to our area of expertise, we only consider articles that use ML to output an approximate solution to one or more fluid-related PDEs.
The PDEs we include in the review are: advection, advection-diffusion, Burgers', Euler, Navier-Stokes, reacting Navier-stokes, advection-diffusion-reaction, Korteweg–de Vries (KdV), Kuramoto–Sivashinsky (KS), shallow water, parabolized stability equations, Poisson, wave, elliptic diffusion, Helmholtz, Laplace, convective heat transfer, plasma models including MHD, PIC, \& Hasegawa-Wakatani, particle-based fluid dynamics, Boltzmann or plasma collision operators, \& black hole hydrodynamics. 

We only include articles that compare the speed, computational cost, or some proxy for speed, of an ML-based solver to that of a standard numerical method used to solve that PDE. Examples of proxies for speed or cost include number of iterations or resolution in space or time. The comparison must be made in a figure, in a table, or in a quantitative statement in the text. A qualitative statement (e.g., ``our method is more efficient'') counts as a valid comparison if it is supported by quantitative or visual supporting evidence.

We define the `primary outcome' as follows. First, if the article has a quantitative comparison (e.g., ``$56\times$ faster'' or ``$4\times$ coarser'') in the abstract, we use that comparison. If no quantitative comparisons are made in the abstract, we look for quantitative comparisons in the conclusion, followed by the introduction, followed by the main text. If no quantitative comparison is made in the entire text, then we look for a qualitative comparison (e.g., ``significantly faster'' or ``outperforms'') beginning in the abstract, followed by the conclusion, followed by the introduction. If there are multiple quantitative or qualitative comparisons, we use our best judgement to determine which should count as the `primary outcome.'

We ended the search process on April 1st, 2023 and thus only include articles available on or before that date. We didn't find any articles published before 2016 that matched our inclusion criteria. Tables \ref{tab:baselines} and \ref{tab:reproducebaselines} show the number of citations each article has according to Google Scholar as of July 3rd, 2023.

%Abbreviations in table \ref{tab:baselines} are for the following PDEs: incompressible Navier-Stokes (INS), compressible Navier-Stokes (CNS), Navier Stokes Airfoil (NSA), Navier-Stokes cylindrical flow (NSC), Burgers' (BG), wave (WV), advection (AV), advection-diffusion (AVD), Korteweg–De Vries (KdV), Kuramoto–Sivashinsky (KS), compressible Euler (EL), Poisson (PS), parabolized stability equations (PSE), plasma particle-in-cell (PIC), reduced magnetohydrodynamics (RMHD), non-linear poisson (NPS), electroconvection (EC), Boltzmann collision operator (BCO), black hole hydrodynamics (BH), Helmholtz (HH), convection (CV), Laplace (LP). 

\subsubsection{Exclusion criteria}

We exclude from the review articles that only consider PDEs related to solid mechanics, quantum mechanics, multiscale modeling, or other non-fluid-related topics. We exclude Reynolds-averaged Navier-Stokes (RANS) and large-eddy simulation (LES).
We also exclude the following PDEs and problems: weather, climate, Schrodinger, fractional, multiphase flows including gas-particle flow, Darcy flow, reaction-diffusion, Eikonal, parabolic diffusion, very high-dimensional PDEs, Compton scattering, meta-materials, hyper-elasticity, ice flow, vessel dynamics, and CO2 injection. 
We also exclude review articles, theses, presentations, technical reports, articles published in languages other than English, ill-posed \& inverse problems, backstepping \& control problems, surrogates for macroscopic quantities, stochastic differential equations, and ordinary differential equations (ODEs).
We exclude model order reduction (MOR) methods, including SVD-based methods such as proper orthogonal decomposition (POD), Sparse Identification of Non-linear Dynamics (SINDy) and Dynamic Mode Decomposition (DMD). 
We don't exclude kernel-based methods, though we didn't find any kernel-based solvers matching our inclusion criteria. 
We exclude physics-informed neural network (PINN)-based methods, because (a) standard numerical methods are known to outperform PINNs for solving forward problems \cite{grossmann2023can}, (b) the PINN literature is too vast to comprehensively review (e.g., \citet{raissi2019physics} has over 9,000 citations), and (c) we only know of a few articles that have ever reported superior performance with PINNs compared to standard numerical methods, and to the best of our knowledge all of these articles either compare to weak baselines or fail to account for the PINN optimization time.
Of the articles included in our systematic review, it turns out that ``machine learning'' invariably involved the use of neural networks and/or deep learning. 

We exclude articles that compare to no baselines or articles that compare to ML baselines but not a standard numerical method as a baseline.
We exclude articles that compare the accuracy of ML-based solvers with standard numerical methods but not the speed or computational cost.
We exclude six articles \cite{ray2018artificial,wang2023long,ovadia2021beyond,li2021learning,ni2023numerical,dong2021modified} that make a qualitative statement of comparison (e.g., ``our method is more efficient'') that are not supported by quantitative or visual supporting evidence about the relative computational cost.
We exclude five articles \cite{mueller2022leveraging,wang2021fast,schwander2021controlling,donon2020deep,wan2023evolve} that might implicitly be suggesting that their proposed method is more efficient on a fluid-related PDE but never make an explicit statement or comparison about the relative speed (or a proxy for speed).
We excluded four articles \cite{di2023neural,kovachki2021neural,holl2020phiflow,nemmen2020first,wandel2020unsupervised} for having duplicate results with other articles.
We excluded one article \cite{haridas2022deep} that uses neural networks to correct floating-point errors in a 16-bit simulation.

\subsubsection{Search process}

The process of systematically searching for every article matching our search criteria was informal at first, but eventually turned into a formal process that happened in two stages. In the first stage, we compiled in list A the names of every author we knew of who worked on ML and PDEs. For each author in list A, we used their Google Scholar profile to look at every title of every article published since 2016. If the title seemed potentially relevant to ML and PDE solving, we read the abstract. If the abstract suggested that the article might possibly satisfy our search criteria, we added it to list B. 
In the second stage, we used Google Scholar to find every article that cites one of two key articles \cite{kochkov2021machine,brandstetter2022message}. We read the title and abstract of each article. If the abstract suggested that the article might satisfy our search criteria, we added it to list B. 
We also tried using Google Scholar to search for key words such as ``machine learning physics", ``machine learning partial differential equations'', ``machine learning fluids'', ``machine learning accelerate pde'', etc. This third approach did not discovery any new articles that were not already added to list B.

For every article added to list B, we read the introduction and conclusion. We also searched the text for key words such as ``fast'', ``speed'', ``improve'', ``pde'', ``equation'', ``compare'', etc., to determine whether the article might have matched the inclusion criteria. Articles that were once again deemed to potentially match the inclusion criteria were read fully to determine whether they should be included in the systematic review. If we found a citation to a new article that might match the inclusion criteria, we added it to list B as well. We also added every co-author of every article that matched the inclusion criteria to list A.

We didn't count the exact number of titles or abstracts we read in total. We added 258 authors to list A and 358 articles to list B. 82 of the articles in list B matched our inclusion criteria.

While we did our best to find every article matching our search criteria, it is possible we missed some articles. If we missed any articles, they are likely articles with fewer citations and/or articles that didn't cite a few key articles.

\subsubsection{Criteria for evaluating baselines \label{sec:methods_evaluation_rules}}

We introduce two necessary but not sufficient conditions (rules 1 and 2) which must be satisfied to ensure a fair comparison between a ML-based PDE solver and a standard numerical method. Rule 1 is to make comparisons at either equal accuracy or equal runtime. Rule 2 is to compare to an efficient numerical method. These rules are discussed in detail in section \ref{sec:weakbaselines}.

We also introduce three recommendations that we recommend following, but do not require that they be satisfied to ensure a fair comparison.
Recommendation 1 is to be cautious of comparing between general-purpose tools and specialized algorithms.
In order to solve a wide class of PDEs, general-purpose libraries are forced to make choices that trade off efficiency for robustness, making them suboptimal for many PDEs.
In contrast, ML-based solvers are specialized algorithms trained to be optimal for a specific PDE or a narrow class of problems. 
Comparisons between a specialized ML-based solver and a general-purpose library are likely to be unfairly biased in favor of the specialized solver.
None of the articles we found explicitly mentioned any reasons to be wary of comparisons between general-purpose tools and the highly specialized ML-based solver.
Moving forward, we encourage articles to be cautious about making these comparisons and to warn readers of the potential for an unfair comparison.

Recommendation 2 is to justify why the choice of hardware (CPU/GPU/TPU) used for comparison is fair.
Some methods, including neural networks, achieve significant reductions in runtime using graphics processing units (GPU) or tensor processing units (TPU) rather than central processing units (CPU).
Other methods achieve only minimal speedups, or no speedup at all, using GPU/TPU compared to CPU. 
Some methods are not implemented on GPU/TPU, only CPU.
In practice, what type of hardware to use for a fair comparison can be context-dependent and to some extent subjective.
Usually, GPU-GPU or TPU-TPU comparisons will be most fair, but in some contexts CPU-CPU or CPU-GPU comparisons can be considered fair.
Most of the articles we found made reasonable choices for the hardware used when comparing different methods.
Moving forward, rather than making definitive rules regarding the choice of hardware, we encourage articles to explain why they chose the hardware they did and to justify why that choice is fair.

In order to account for the cost of generating data and training models, recommendation 3 is to report the number of surrogate evaluations $N$ needed to reduce the total computational cost in downstream applications. $N$ is only defined if the ML-based solver is faster than the numerical baseline. $N$ can be computed using the formula
\begin{equation}
    C_{\textnormal{data}} + C_{\textnormal{train}} + N \frac{t_B}{s} = N t_B
\end{equation}
where $C_{\textnormal{data}}$ is the time required to generate the training data, $C_{\textnormal{train}}$ is the time required to train the model(s), $t_B$ is the time required for the standard numerical method baseline to compute one surrogate evaluation at equivalent accuracy to the ML model, and $s$ is the speedup of the ML-based solver relative to the numerical baseline.

We now explain which standard numerical methods we consider state-of-the-art for elliptic and advection-dominated PDEs.
For elliptic PDEs (such as Poisson’s, Laplace’s, or Helmholtz equations) finite element methods (FEM) are standard; direct solvers such as LU decomposition are most efficient for small problems, while iterative solvers are more efficient for large problems. For elliptic PDEs, multigrid solvers are typically state-of-the-art for large problems. We suggest using Eigen \cite{guennebaud2010eigen} for LU decomposition and HYPRE \cite{falgout2002hypre} for multigrid methods, though other libraries can also be extremely efficient. 
For advection-dominated PDEs, we recommend using second or higher-order shock-capturing finite volume (FV) methods for problems with shocks (such as compressible Navier-Stokes and Burgers’ equations), while using higher-order methods for problems with smooth solutions (such as the advection equation with smooth solutions, the incompressible Navier-Stokes equations or the compressible Navier-Stokes equations in the weakly compressible limit). For advection-dominated problems with smooth solutions, pseudo-spectral methods are usually state-of-the-art when applicable, though discontinuous Galerkin (DG) methods are also extremely efficient. We have found that higher-order DG methods (polynomial order 2 or higher) work better than lower-order DG methods (order 0 or 1), though there are diminishing returns for using very high-order (polynomial order 3 or higher) DG methods. Moving forward, for the Navier-Stokes equations we recommend comparing the performance of FV, DG, and (if applicable) spectral methods and choosing the strongest baseline.
In general, first-order methods should not be used as baselines for fluid-related PDEs. First-order methods tend to be extremely diffusive and require high grid resolution.
Explicit time-stepping schemes are usually preferred for advection-dominated PDEs; for these problems, using implicit time-stepping typically leads to inefficient numerical methods. The time-step restriction from the Courant–Friedrichs–Lewy (CFL) condition is usually sufficiently small that the dominant error is from spatial discretization, and so the choice of explicit time-stepping scheme is less important than the spatial discretization.
However, in some cases (such as fluid-structure interaction \cite{mayr2018adaptive}) the temporal discretization errors can dominate, in which case specialized time-stepping schemes that compute error estimates and perform adaptive time-stepping perform best \cite{reynolds2014arkode}.

To help ensure that we applied these rules fairly, we emailed the authors of each article to give them an opportunity to point out any errors we might have made in applying rules 1 and 2. We received 15 responses about 23 articles. Seven responses expressed agreement and gave suggestions for improvement, two provided additional information, and six expressed disagreement. Based on the responses, we modified six entries in Table \ref{tab:baselines} and one entry in Table \ref{tab:reproducebaselines}.

\subsubsection{Details of stronger baselines in Table \ref{tab:reproducebaselines}}

For articles 1, 4, and 12 we use a Runge-Kutta discontinuous galerkin (RKDG) method to solve the 2D incompressible Navier-Stokes equations on a periodic domain. We use a third-order strong-stability preserving (SSPRK3) ODE integration \cite{gottlieb2001strong}. We use second-order discontinuous polynomial basis functions with a serendipity basis function which uses 8 basis functions per grid cell. We use LU decomposition and a continuous Galerkin formulation to solve the Poisson equation at each Runge-Kutta stage. The full scheme is explained in \citet{hakim2019discontinuous}. \citet{dresdner2022learning} also reproduces article 4 using a pseudo-spectral implementation. The details can be found at \url{https://github.com/google/jax-cfd}.
For articles 2 and 6 we use again a RKDG scheme, except this time to solve the 1D advection and 1D Burgers' equations with periodic boundary conditions. We again use SSPRK3 ODE integration and second-order Legendre polynomial basis functions. The full scheme is explained in \citet{cockburn1989tvb}. For article 43 we again use an RKDG scheme to solve the 1D Burgers' equation, except with Dirichlet instead of periodic boundary conditions.
For articles 8 and 15 we use a finite volume scheme with Godunov flux \cite{durran2013numerical} and SSPRK3 ODE integration to solve the 1D Burgers' equation with periodic boundary conditions.
For articles 14 and 68 we solve the Poisson equation on a square periodic domain using a continuous Galerkin formulation with linear basis functions. We use an LU decomposition to solve the resulting linear system.
For article 15 we solve the 1D wave equation with Dirichlet boundary conditions, using SSPRK3 timestepping and a finite-volume method with irregular grid spacing.

\subsection{Random sample of ML-for-PDE articles}

To approximate a random sample of articles in ML-for-PDE solving, we use Google Scholar to find all 400 articles that cite \citet{kochkov2021machine} as of June 18th, 2023. We include articles whose abstracts mention positive and/or negative experimental results. We define a result as: proposing a method to tackle a problem and using quantitative metrics or qualitative descriptions to describe the performance of the method on the problem.
We classify each abstract based on whether it mentions positive and/or negative results.
Many articles have the pattern `method A has negative aspects, we introduce method B which solves those negative aspects'; we don't consider this pattern a negative result unless the article discusses negative results associated with method B. 
Three articles comment that a method has limitations or is limited in scope; we do not count those comments as negative results. We exclude review articles, theses, duplicates, articles written in languages other than English, articles that don't mention experimental results anywhere in the abstract, and articles that are unrelated to ML or statistical learning. 
We excluded ten articles which mention experimental results, but don't give any indication as to whether those results are positive or negative.
232 out of 400 articles (58\%) in the sample match our inclusion criteria.

\subsection{Random sample of PINN articles}

Physics-informed neural networks (PINNs) \cite{raissi2019physics} are a popular method which can solve PDEs and inverse problems associated with PDEs. We look for a random sample of articles which (a) use PINNs to solve a fluid-related PDE, (b) which focus on solving ``forward'' problems and not inverse or ill-posed problems, (c) which report the accuracy of the PINN-based solver, and (d) which generate their own data using a standard numerical method or general-purpose solver to measure accuracy. We define ``fluid-related'' using the same inclusion and exclusion criteria defined earlier. Once again, we exclude review articles, theses, duplicates, technical reports, and articles not written in English. 
We also exclude articles that use an analytic solution to measure accuracy.
If an article uses PINNs to solve both forward and inverse (i.e., ill-posed) problems, we only include the article if a majority of the problems it solves are forward problems.

To obtain a random sample of articles matching these criteria, we use Google Scholar to search within all 5,640 articles which cite \citet{raissi2019physics} as of June 21st, 2023. Using the ``search within citing articles'' function of Google Scholar, we type into the search bar

\begin{quote}
    ``PINN" AND ``Burgers" OR ``navier" OR ``stokes" OR ``fluid" OR ``advection" OR ``KdV" OR ``Kuramoto" OR ``Sivashinsky" OR ``Euler" OR ``convection" OR ``Laplace" OR ``poisson" OR ``parabolized stability" OR ``plasma" OR ``collision" OR ``MHD" OR ``Helmholtz"
\end{quote}
and sort by relevance. This returns 1,000 articles, which is the most articles Google Scholar will return in a single search. We focus only on the first 250 articles, sorting by relevance rather than by date. If the title and abstract potentially matches the inclusion criteria, we add it to a list of potentially relevant articles. We added 155 articles (62\%) to the list of potentially relevant articles. We then read each article in the list closely to see if it matches our inclusion criteria. 
We exclude eight articles which claim to outperform standard numerical methods with PINNs, though each of these articles compares to a weak baseline or, more often, doesn't account for the PINN optimization time.
This search process ultimately returns 60 articles, which we use as sample B.

It is well known that PINN-based solvers are slower than standard numerical methods (except possibly when meta-learning is used, e.g., in \citet{qin2022meta}) \cite{grossmann2023can,karniadakis2021physics}.
Importantly, we can assume that the authors of each article in sample B know this. We can assume this because (a) is it common knowledge, (b) we only include articles that generate their own data from standard numerical methods to measure accuracy, so they must have known the runtime of that method, and (c) most of the articles in sample B are focused on developing strategies to improve the speed of PINNs, and thus these authors recognize that speed is a limitation of PINNs. Using PINNs to get good accuracy on most PDEs takes (in most cases) hours to days, while doing so with a standard numerical method takes between a fraction of a second to minutes \cite{chuang2022experience,chuang2023predictive,grossmann2023can}, depending on the PDE and the solver used. The authors doubtlessly noticed the difference in runtime between the PINN and the standard solver, and must have known (at least implicitly) that the PINN-based solver was less efficient than the standard numerical method.

\backmatter

\section{Declarations}

\bmhead{Data availability}

The lists of authors and articles generated during the systematic review and the categorizations of every article in the random samples are publicly available at \url{https://doi.org/10.17605/OSF.IO/GQ5B3} \cite{McGreivy_dataset_2024}.

\bmhead{Code availability}

The code required to reproduce the results in table \ref{tab:reproducebaselines} is available at \url{https://github.com/nickmcgreivy/WeakBaselinesMLPDE/} \cite{McGreivy_Github_2024}. We provide instructions for running the code and interpreting the results.

\bmhead{Authors' contributions}

N.M. conceptualized the systematic review, searched for papers matching the inclusion criteria, evaluated rule 1 and rule 2 for each article, conceptualized and carried out analyses to measure the effect of reporting biases, and wrote the code and the manuscript. A.H. designed strong baselines, evaluated rule 2 for each PDE, provided instructions for implementing the code, edited the manuscript, and supervised the research.

\bmhead{Supplementary information} Supplementary Information is available for this paper. 

\bmhead{Acknowledgments}
N.M. was supported via DOE contract DE-AC02-09CH11466 for the Princeton Plasma Physics Laboratory. A.H. was supported by the Partnership for Multiscale Gyrokinetic Turbulence (MGK) and the High-Fidelity Boundary Plasma Simulation (HBPS) projects, part of the U.S. Department of Energy (DOE) Scientific Discovery Through Advanced Computing (SciDAC) program, and the DOE’s ARPA-E BETHE program, via DOE contract DE-AC02-09CH11466 for the Princeton Plasma Physics Laboratory. 
The funders had no role in study design, data collection and analysis, decision to publish or preparation of the manuscript.

\bmhead{Competing interests} 

The authors declare no competing interests.

\newpage

\appendix

\section*{Supplementary Information}

We now give detailed explanations for each of the 76 articles in table 1, as well as the 6 articles that claim to either underperform or have varied performance relative to a standard numerical method.

\appendixentry{\fnoli}{\citefnoli}{Fourier Neural Operator for Parametric Partial Differential Equations}{\fnolicitations}
\appendixdata{1D Burgers', 2D incompressible Navier-Stokes}
{``up to three orders of magnitude faster compared to traditional PDE solvers.''}
{No}
{Pseudo-spectral method}
{{\xmark} ``All data are generated on a $256\times256$ grid and are downsampled to $64\times64$.'' ``On a $256\times 256$ grid, the Fourier neural operator has an inference time of only 0.005s compared to the 2.2s of the pseudo-spectral method used to solve Navier-Stokes.'' This comparison is not at equal accuracy. To compare at equal accuracy, reduce the resolution of the pseudo-spectral method until the two methods have equal accuracy (as measured by table 1).  }
{{\cmark} A pseudo-spectral method or DG method is considered state-of-the-art for the 2D incompressible Navier-Stokes equations with periodic boundary conditions.}
{N/A}
{{\xmark} We replicated the primary outcome of this article using a DG method and found that the speedup of Fourier neural operator on GPU was $7\times$ faster than our laptop CPU, not three orders of magnitude faster.}

\appendixentry{\deeponetlulu}{\citedeeponetlulu}
{Learning nonlinear operators via DeepONet based on the universal approximation theorem of operators}{\deeponetlulucitations}
\appendixdata{1D advection (4 cases), 1D advection-diffusion}
{``as we show in Supplementary Table 5, the computational cost of running inference of DeepONet is substantially lower than for the numerical solver.''}
{No}
{``The reference solutions of all deterministic PDEs are obtained by a second-order finite difference method.'' ``To generate the training dataset, we solve the system using a finite difference method on a 100 by 100 grid.''}
{{\xmark} Table 5 does not list the runtime at equal accuracy, and the grid resolution of the finite difference method is not reduced to match the accuracy of the DeepONet.}
{{\cmark} A second-order finite-difference method is fairly efficient for the 1D advection equation, but note that DG methods are likely more efficient for the 1D advection equation with smooth solutions.}
{N/A}
{{\xmark} We consider the 1D advection equation (case 1). As listed in table S8, this is the linear advection equation from $x \in [0, 1]$ and $t \in [0, 1]$. The error of the $100 \times 100$ second-order finite-difference method is shown in figure S11; the runtime is $9 \times 10^{-3}$ seconds using 1 core on CPU. We instead run a DG code with quadratic basis functions with 13 grid cells with a similar initial condition, the runtime is $6 \times 10^{-5}$ seconds on my laptop which has 2 cores. In summary, we were able to achieve similar accuracy with an order of magnitude lower runtime. }

\appendixentry{\tompson}{\citetompson}
{Accelerating Eulerian Fluid Simulation With Convolutional Networks}{\tompsoncitations}
\appendixdata{2D and 3D Poisson for computer graphics (GPU-only) in real-time (low-accuracy)}
{``For Jacobi to match the divergence performance of our network, it requires 116 iterations and so is $4.1\times$ slower than our network.'' ``Note that for fair quantitative comparison of output residual, we choose the number of Jacobi iterations (34) to match the FPROP time of our network (i.e. to compare divergence at fixed compute).'' Supported by figure 5. }
{No}
{``A Jacobi-based iterative solver and a PCG-based solver (with incomplete Cholesky L0 preconditioner).'' The PCG baseline ``is orders of magnitude slower and has been omitted for clarity.'' }
{{\cmark} Compares divergence (accuracy) at constant compute (speed), and speed at constant accuracy.}
{{\cmark} For Poisson's equation, Multigrid methods (or preconditioners) are highly efficient at high accuracy, because they converge in many fewer iterations than Jacobi. See, e.g., figure 4e of article {\tang} or figure 23 of article \cheng. However, for the scenario considered in the paper of GPU-only real-time computer graphics applications, the authors argue that multigrid methods are less efficient than Jacobi iteration. This argument seems plausible, but ultimately we are unsure if it is correct. }
{N/A}
{{\cmark} }

\appendixentry{\mlacceleratedcfd}{\citemlacceleratedcfd}
{Machine learning–accelerated computational fluid dynamics}{\mlacceleratedcfdcitations}
\appendixdata{2D incompressible Navier-Stokes}
{``our  results  are  as  accurate  as  baseline  solvers  with 8 to 10$\times$ finer resolution in each spatial dimension, resulting in 40- to 80-fold computational speedups.'' Supported by figures 2a and 2b.}
{No}
{Finite-volume method based on a Van-Leer flux limiter.}
{{\cmark} See figure 2b.}
{{\xmark} Pseudo-spectral and DG methods are highly efficient for this problem. The original authors of this article replicated this result using a strong spectral baseline on TPU (see article \dresdner) and found that a pseudo-spectral baseline was much faster than the weaker FV baseline and faster than a similar ML-based solver. On GPU, we find that the PS baseline is over 80x faster than the FV baseline. We also replicated this result using a DG baseline and found that DG-based methods could solve these equations at 10 to 11$\times$ coarser resolution and (on CPU) 20-40$\times$ faster than the original baseline.}
{N/A}
{{\xmark}}

\appendixentry{\meshbasedpfaff}{\citemeshbasedpfaff}
{Learning Mesh-Based Simulation with Graph Networks}{\meshbasedpfaffcitations}
\appendixdata{2D incompressible Navier-Stokes cylindrical flow, 2D compressible Navier-Stokes airfoil wing}
{``Our method is also highly efficient, running 1-2 orders of magnitude faster than the simulation on which it is trained.'' Supported by table 1 and section A.5.1.}
{No}
{COMSOL for incompressible Navier-Stokes, SU2 for compressible Navier-Stokes. Also compares to ANSYS.}
{{\xmark} Table 1 compares the runtime between the (highly accurate) ground-truth solver and the less accurate ML-based solver. This comparison is not at equal accuracy.}
{{\cmark} Though we don't have reason to believe COMSOL and SU2 are inefficient general-purpose tools, we recommend being cautious when evaluating this comparison (see recommendation 1 in Methods).}
{N/A}
{{\xmark}}

\appendixentry{\barsinai}{\citebarsinai}
{Learning data-driven discretizations for partial differential equations}{\barsinaicitations}
\appendixdata{1D Burgers', 1D Korteweg-de Vries, 1D Kuramoto-Sivashinsky}
{``The resulting numerical methods are remarkably accurate, allowing us to integrate in time a collection of nonlinear equations in 1 spatial dimension at resolutions $4\times$ to $8\times$ coarser than is possible with standard finite-difference methods.''}
{No}
{A FV method with a fifth-order upwind-biased WENO scheme with Godunov flux}
{{\cmark} See figure 3c.}
{{\xmark} We consider the 1D Burgers' equation (in figure 3). We replicate figure 3c, and compare WENO to DG methods with polynomial orders 2 and 3. We find (see github) that, as with the ML-based solver, DG methods are able to solve the 1D Burgers' equation at resolutions $4\times$ to $8\times$ coarser than WENO. The ML-based solver is still able to solve the 1D Burgers' equation with 2-4$\times$ fewer degrees of freedom compared to the DG-based method. We give an {\xmark} for rule 2 because we were able to replicate the article's primary outcome and achieve significantly improved performance with a stronger baseline.
}
{N/A}
{\xmark}

\appendixentry{\deepfluids}{\citedeepfluids}
{Deep Fluids: A Generative Network for Parameterized Fluid Simulations}{\deepfluidscitations}
\appendixdata{2D \& 3D incompressible Navier-Stokes (smoke/graphics)}
{``Reconstructed velocity fields are generated up to 700$\times$ faster than re-simulating the data with the underlying CPU solver.'' Supported by table 1.}
{No}
{Mantaflow}
{{\xmark} The resolution of the underlying CPU solver is not reduced to match the accuracy of the ML-based solver.}
{{\cmark} Though we don't have reason to believe that Mantaflow is an inefficient general-purpose tool, we recommend being cautious when evaluating this comparison (see recommendation 1 in Methods).}
{N/A}
{\xmark}

\appendixentry{\deeponetwang}{\citedeeponetwang}
{Learning the solution operator of parametric partial differential equations with physics-informed DeepOnets}{\deeponetwangcitations}
\appendixdata{1D Burgers'}
{``up to three orders of magnitude faster compared a conventional PDE solver.'' Supported by figure 11.}
{No}
{Spectral solver (Chebfun)}
{{\xmark} Figure 11 compares the runtime of a highly accurate spectral solver to that of a less accurate ML-based solver.}
{{\cmark} A spectral solver is highly efficient for Burgers' equation, as long as the diffusion coefficient is sufficiently high (so that shocks are not too strong). }
{N/A}
{{\xmark} We replicate the 1D Burgers' setup using a FV method with Godunov flux. Using 100 gridpoints on CPU, this gives an L2 error of 1\% and a runtime of $4.2\times 10^{-4}$ seconds, over an order of magnitude faster than the ML-based solver. On GPU, we can solve 1000 PDEs in $1.2\times10^{-2}$ seconds, again an order of magnitude faster than the MLP in figure 11b.}

\appendixentry{\solverinthe}{\citesolverinthe}
{Solver-in-the-Loop: Learning from Differentiable Physics to Interact with Iterative PDE-Solvers}{\solverinthecitations}
\appendixdata{2D Burgers', 2D Poisson, 2D and 3D incompressible Navier-Stokes wake dynamics}
{``A speed-up of more than 68$\times$ [for the simulation in figure 1].'' Supported by appendix C, section titled ``runtime performance''.}
{No}
{Reference simulation (baseline appears to be the \href{https://github.com/tum-pbs/PhiFlow}{PhiFlow} library, which is based on a MAC grid data structure). See section B.5 for details about 3D setup, and section B.1 for details about 2D baseline. Baseline appears to be a FV method.}
{{\xmark} The 68$\times$ speedup compares the reference simulation with MAE of 0 to a reference simulation with MAE of 0.13. This is only a 28\% improvement over the `source' simulation, which has an MAE of 0.167. A fair comparison would be between two simulations with approximately equal MAE.}
{{\cmark} While we recommend comparing with both FV and DG methods for the Navier-Stokes equations, we consider FV highly efficient for problems with fluid-structure interaction. Though we don't have reason to believe that PhiFlow is an inefficient general-purpose tool, we also recommend being cautious when evaluating this comparison (see recommendation 1 in Methods). }
{N/A}
{\xmark}

\appendixentry{\deepmmnetcvt}{\citedeepmmnetcvt}
{DeepM\&Mnet: Inferring the electroconvection multiphysics fields based on operator approximation by neural networks}{\deepmmnetcvtcitations}
\appendixdata{2D electroconvection (steady state)}
{``The speedup of DeepONets prediction versus the NekTar simulation for forward independent conditions is about 10,000 folds.'' Not supported by any other evidence. }
{No}
{Nektar: high-order spectral element (3 modes), 5 quadrature points in each direction, with second-order stiffly stable timestepping scheme. `` }
{{\xmark} To make a fair comparison, reduce the resolution of the Nektar simulation (below $32\times 32$) until its accuracy is equal to that of the DeepONet.}
{{\cmark} Though we don't have reason to believe that Nektar is an inefficient general-purpose tool, we recommend being cautious when evaluating this comparison (see recommendation 1 in Methods).}
{}
{\xmark}

\appendixentry{\belbuteperez}{\citebelbuteperez}
{Combining Differentiable PDE Solvers and Graph Neural Networks for Fluid Flow Prediction}{\belbuteperezcitations}
\appendixdata{2D compressible Navier-Stokes airfoil wing}
{``the substantial speedup of neural network CFD predictions.'' Supported by Table 1.}
{No}
{SU2}
{{\xmark} Compares ground truth (runtime 137s, RMSE 0.0) to CFD-GCD (runtime 2.0s, RMSE $5.4\times10^{-2}$) instead of comparing at equal accuracy.}
{{\cmark} Though we don't have reason to believe that SU2 an inefficient general-purpose tool, we recommend being cautious when evaluating this comparison (see recommendation 1 in Methods).}
{}
{\xmark}

\appendixentry{\pinoli}{\citepinoli}{Physics-Informed Neural Operator for Learning Partial Differential Equations}{\pinolicitations}
\appendixdata{1D Burgers', 2D incompressible Navier-Stokes}
{``Further, in PINO, we incorporate the Fourier neural operator (FNO) architecture which achieves orders-of-magnitude speedup over numerical solvers." Supported by figure 8.}
{No}
{Same spectral solver as in article \fnoli.}
{{\cmark} See figure 8.}
{{\xmark} The transient flow problem is identical to that in article \fnoli, with one key difference: the Reynolds number is now 20, instead of $10^3$-$10^5$. Thus, the problem is now diffusion-dominated rather than advection-dominated. We use a DG solver with second-order polynomial basis functions to replicate this result, except we reduce the resolution to $3\times3$ and change the timestep accordingly. With a $3\times3$ resolution, we find an error of 2-3\% with a runtime of 0.035s, 7$\times$ slower than the PINO method with similar accuracy.}
{N/A}
{\xmark}

\appendixentry{\projectionyang}{\citeprojectionyang}{Data-driven projection method in fluid simulation}{\projectionyangcitations}
\appendixdata{3D Poisson}
{``Experimental results demonstrated that our data-driven method drastically speeded up the computation in the projection step.'' Supported by table II.  }
{No}
{Preconditioned Conjugate Gradient (PCG) linear solver}
{{\xmark} The ML-based solver doesn't use iteration. Thus, the accuracy of the ML-based solver (as evidenced by figure 4) isn't as high at that of the PCG baseline.}
{{\xmark} For Poisson's equation, Multigrid methods (or preconditioners) are highly efficient. See, e.g., figure 4e of article {\tang} or figure 23 of article \cheng.}
{}
{\xmark}

\appendixentry{\neuralconverge}{\citeneuralconverge}{Learning Neural PDE Solvers with Convergence Guarantees}{\neuralconvergecitations}
\appendixdata{2D Poisson}
{``[Our model] achieves 2-3 times speedup compared to state-of-the-art solvers.'' Supported by figure 2.}
{No}
{``The FEniCS model is set to be the minimal residual method with algebraic multigrid preconditioner, which we measure to be the fastest compared to other methods such as Jacobi or Incomplete LU factorization preconditioner.'' }
{\cmark ``We evaluate the convergence rate by calculating the computation cost required for the error to drop below a certain threshold.'' This article compares runtime at equal accuracy. }
{{\xmark} Consider figure 2a: FEniCS takes almost 20 seconds to solve a $256\times 256$ Poisson problem on a square domain. We implement Poisson's equation on a square domain using a direct solve (LU decomposition) and find that the direct solve takes 12 milliseconds, over three orders of magnitude faster than FEniCS. Multigrid is a weak baseline relative to direct methods for sufficiently small problems.}
{N/A}
{\xmark}

\appendixentry{\messagepassing}{\citemessagepassing\hspace{0.1cm}(Note: we consider version 2 of this article on ArXiv, version 3 was uploaded after private communication with the authors.)}{Message Passing Neural PDE Solvers}{\messagepassingcitations}
\appendixdata{1D Burgers', 1D wave, 2D incompressible Navier-Stokes (smoke/graphics)}
{``Our model outperforms state-of-the-art numerical solvers in the low resolution regime in terms of speed and accuracy.'' Supported by tables 1 and 2.}
{No}
{WENO5 (Burgers) and spectral (wave)}
{{\cmark} See tables 1 and 2.}
{{\xmark} It's a slow implementation of WENO5 and spectral. }
{N/A}
{{\xmark} We replicated both problems with a stronger baseline, and found that the stronger baseline was orders of magnitude faster than the ML-based solver.}

\appendixentry{\deepmmnet}{\citedeepmmnet}{DeepM\&Mnet for hypersonics: Predicting the coupled flow and finite-rate chemistry behind a normal shock using neural-network approximation of operators}{\deepmmnetcitations}
\appendixdata{2D reacting Navier-Stokes}
{``DeepONets can be over five orders of magnitude faster than the CFD solver employed to generate the training data.'' }
{No}
{CFD solver (no details given) coupled with the MUTATION library}
{{\xmark} This article compares a method with relative MSE of $1e-5$ to a CFD solver with relative MSE of 0.0. For a fair comparison, reduce the resolution of the CFD solver until the relative MSE is equal. }
{{\xmark} ``Even though variations only take place along the streamwise direction, the actual computations were performed in two dimensions as stated earlier.'' The two dimensions are streamwise ($x=x_1$) and normal ($y=x_2$). In other words, this article uses a 2D code as a baseline for a 1D problem. }
{}
{\xmark}

\appendixentry{\frameworkmishra}{\citeframeworkmishra}{A machine learning framework for data driven acceleration of computations of differential equations}{\frameworkmishracitations}
\appendixdata{1D Burgers', advection, Euler}
{``Numerical experiments involving both linear and non-linear ODE and PDE model problems demonstrate a significant gain in computational efficiency over standard numerical methods.'' Supported by tables 5, 6, and 7, as well as page 20.}
{No}
{Rusanov for Burgers' and Euler, backwards euler time-stepping for advection}
{{\cmark} Focuses on speedup at constant error (page 20) or error at constant speed (see ``gain'' on page 12).}
{{\xmark} The Rusanov scheme is a first-order scheme, and is much more diffusive and less accurate at solving the 1D Euler equations than a higher-order scheme (e.g., MUSCL scheme with reconstruction in characteristic variables). For Burgers' equation, WENO5 would be a strong FV baseline.}
{N/A}
{\xmark}

\appendixentry{\optimizemultigrid}{\citeoptimizemultigrid}{Learning to Optimize Multigrid PDE Solvers}{\optimizemultigridcitations}
\appendixdata{2D elliptic diffusion}
{``Experiments on a broad class of 2D diffusion problems demonstrate improved convergence rates compared to the widely used Black-Box multigrid scheme.'' Supported by figure 3.}
{No}
{Black box multigrid scheme}
{{\cmark} Measures number of iterations (i.e., convergence rate) at constant accuracy. See figure 3.}
{{\cmark} Multigrid is highly efficient for elliptic problems.}
{N/A}
{\cmark }

\appendixentry{\donga}{\citedonga}{Local extreme learning machines and domain decomposition for solving linear and nonlinear partial differential equations}{\dongacitations}
\appendixdata{1D Helmholtz, 1D advection, 1D diffusion, 1D non-linear helmholtz, 1D Burgers', 2D Poisson.}
{``The computational performance of the current method is on par with, and often exceeds, the FEM performance in terms of the accuracy and computational cost." Supported by tables 2, 5, 7, 10, 11, as well as figures 50 and 52.}
{No}
{FEniCS, linear elements, second-order, sparse LU decomposition for linear solver, Newtons method for non-linear equations. For 1D diffusion and 1D Burgers', uses second-order backward differentiation formula (BDF2), an implicit timestepping method. }
{{\cmark} Figures 50 and 52 make plots of speed versus accuracy, allowing readers to make comparisons at equal accuracy or equal speed.}
{{\xmark} The primary outcome (in this case defined by what is reported in the abstract) compares the ML-based solver to a weaker baseline (second-order FEM). In the appendix, a stronger baseline is compared to (higher-order FEM). Both the ML-based solver and second-order FEM underperform relative to higher-order FEM. Because the primary outcome reports performance relative to the weaker baseline rather than the stronger baseline, we determine that rule 2 is not satisfied. To satisfy rule 2, report in the abstract performance relative to the strong baseline.}
{N/A}
{\xmark}

\appendixentry{\rayb}{\citerayb}{Detecting troubled-cells on two-dimensional unstructured grids using a neural network}{\raybcitations}
\appendixdata{2D advection, Burgers', Euler}
{``Through several numerical tests, the MLP indicator has been shown to outperform the TVB indicator (with various values of M) both in terms of solution accuracy and the number of cells flagged, while maintaining a comparable computational cost.'' Supported by table 6.}
{No}
{TVB limiter}
{{\cmark} The computational time in table 6 is nearly constant. Comparing accuracy at constant runtime.}
{{\cmark} TVB limiter is highly efficient.}
{}
{\cmark }

\appendixentry{\novelcnn}{\citenovelcnn}{A Novel CNN-Based Poisson Solver for Fluid Simulation}{\novelcnncitations}
\appendixdata{3D Poisson (smoke)}
{``We have shown that our approach accelerates the projection step in a conventional Eulerian fluid simulator by two orders of magnitude.'' Supported by table 1 (page 1462) and section 7.2.}
{No}
{Preconditioned conjugate gradient, both with and without multigrid preconditioning.}
{{\xmark} Compares runtime at different accuracy. ``The accuracy of the PCG-based method [is] higher than our CNN-based solver. As such, our CNN-based linear equation solver is suitable for the simulations which are not strict with numerical accuracy.''}
{{\xmark} The primary outcome of this article (in this case defined by what is reported the abstract) compares to the weaker baseline (MIC(0)-PCG), while table 1 compares to both weak and strong baselines (MG-PCG). Because the primary outcome reports performance relative to the weaker baseline rather than the stronger baseline, we determine that rule 2 is not satisfied. To satisfy rule 2, report in the abstract performance relative to the strong baseline.}
{N/A}
{\xmark}

\appendixentry{\wandel}{\citewandel}{Learning Incompressible Fluid Dynamics from Scratch -- Towards Fast, Differentiable Fluid Models that Generalize}{\wandelcitations}
\appendixdata{2D incompressible Navier-Stokes wake dynamics}
{``The $\vec{v}$-Net as well as the $\vec{a}$-Net'' are significantly faster than PhiFlow ($11\times$ on CPU and $40\times$ on GPU)." Supported by table 1 and Appendix E.}
{No}
{Phiflow (FV method, relies on iterative conjugate gradient solver) using $100\times 100$ grid.}
{{\xmark} The loss (which in this case, doesn't measure accuracy because the pressure is evolved independently of the velocity when pressure is not independent of velocity) is compared, but even so, the runtime is not compared at equal values of the loss.}
{{\cmark} While we recommend comparing with both FV and DG methods for the Navier-Stokes equations, we consider FV highly efficient for problems with fluid-structure interaction. Though we don't have reason to believe that PhiFlow is an inefficient general-purpose tool, we also recommend being cautious when evaluating this comparison (see recommendation 1 in Methods).}
{N/A}
{{\xmark}}

\appendixentry{\shan}{\citeshan}{Study on a Fast Solver for Poisson’s Equation Based on Deep Learning Technique}{\shancitations}
\appendixdata{2D \& 3D Poisson}
{``Numerical experiments show that the same ConvNet architecture is effective for both 2-D and 3-D models\dots with a significant reduction in computation time compared to the finite-difference solver.'' Supported by last paragraph before conclusion.}
{No}
{Finite-difference method }
{{\xmark} The convolutional network presumably has lower accuracy than the finite-difference baseline, because the authors never mention the number of iterations of the finite-difference baseline. The authors admit that their ``comparison is not exactly fair.'' }
{{\xmark} For Poisson's equation, Multigrid methods (or preconditioners) are highly efficient. See, e.g., figure 4e of article {\tang} or figure 23 of article \cheng. Note also that for sufficiently small 2D problems, direct solves (such as LU decomposition) will outperform Multigrid methods. Although this article says nothing about how the system of linear equations (equation 7) is solved, the baseline takes 17s to solve the 2D problem on a $64\times64$ grid, which is orders of magnitude slower a direct solve (such as LU decomposition).}
{N/A}
{{\xmark }}

\appendixentry{\algebraic}{\citealgebraic}{Learning Algebraic Multigrid Using Graph Neural Networks}{\algebraiccitations}
\appendixdata{2D elliptic diffusion}
{``Experiments on a broad class of problems demonstrate improved convergence rates compared to classical AMG.'' Supported by table 4.}
{No}
{algebraic multigrid}
{{\cmark} Table 4: ``required to reach specified tolerance.'' Compares number of iterations at constant accuracy.}
{{\cmark} Algebraic multigrid is highly efficient for elliptic problems.}
{N/A}
{\cmark}

\appendixentry{\zhuang}{\citezhuang}{Learned discretizations for passive scalar advection in a 2-D turbulent flow}{\zhuangcitations}
\appendixdata{2D advection}
{``The method maintains the same accuracy as traditional high-order flux-limited advection solvers, while using 4× lower grid resolution in each dimension.''}
{No}
{second-order Van Leer advection scheme (FV)}
{{\cmark} See figure 8.}
{{\xmark} DG or pseudo-spectral methods are state-of-the-art for scalar advection. DG schemes in particular will solve the advection equations at much coarser resolution. }
{N/A}
{\xmark}

\appendixentry{\pathak}{\citepathak}{Using Machine Learning to Augment Coarse-Grid Computational Fluid Dynamics Simulations}{\pathakcitations}
\appendixdata{2D incompressible Navier-Stokes}
{``The ML-assisted coarse-grid evolution resulted in corrected solution trajectories that were consistent with the solutions computed at a much higher resolution in space and time.'' Supported by figures 2 and 3.}
{No}
{Dedalus, a spectral solver.}
{{\cmark} Although this article never directly reduces the runtime to make a comparison at equal accuracy, it is fair to say that the solver runs at `lower resolution.' If this article had said `$4\times$ lower resolution', that would not have been fair.}
{{\cmark} Spectral solvers are highly efficient for the incompressible Navier-Stokes equations.}
{{\cmark} Though we don't have reason to believe that Dedalus is an inefficient general-purpose tool, we recommend being cautious when evaluating this comparison (see recommendation 1 in Methods).}
{\cmark}

\appendixentry{\leoni}{\citeleoni}{DeepONet prediction of linear instability waves in high-speed boundary layers}{\leonicitations}
\appendixdata{Parabolized stability equations}
{``\dots at a very small computational cost compared to discretization of the original equations.'' Supported by subsection V.B.}
{No}
{The code in reference [5]}
{{\xmark} Section V.B compares the runtime of the highly accurate code (error 0.0) to the less accurate ML-based forward solver (error 2-5\%, see figure 10). This comparison is not at equal accuracy.}
{{\cmark} We are unsure if this code is an efficient numerical method.}
{N/A}
{\xmark}

\appendixentry{\fnodeformation}{\citefnodeformation}{Fourier Neural Operator with Learned Deformations for PDEs on General Geometries}{\fnodeformationcitations}{}
\appendixdata{2D compressible Navier-Stokes airfoil wing, 2D incompressible Navier-Stokes pipe flow}
{``Geo-FNO is $10^5$ times faster than the standard numerical solvers.'' }
{No}
{Second-order implicit FV solver.}
{{\xmark} The comparison is made between the highly accurate (but slow) ground truth solver and the less accurate ML-based solver.}
{{\cmark} Finite-volume methods are highly efficient for the 2D compressible Navier-Stokes airfoil problem.}
{N/A}
{\xmark}

\appendixentry{\stevensa}{\citestevensa}{Enhancement of shock-capturing methods via machine learning}{\stevensacitations}
\appendixdata{1D Burgers', advection, Euler}
{``We find that our method outperforms WENO in simulations where the numerical solution becomes overly diffused due to numerical viscosity.'' Supported by figure 8.}
{No}
{WENO}
{{\cmark} Faster runtime at equal accuracy (see figure 8).}
{{\cmark} WENO5 is highly efficient for problems with shocks/discontinuities.}
{N/A}
{\cmark }

\appendixentry{\illarramendi}{\citeillarramendi}{Towards an hybrid computational strategy based on Deep Learning for incompressible flows}{\illarramendicitations}
\appendixdata{2D Poisson}
{``For the same accuracy to be achieved with only Jacobi iterations, the calculation is 3.2 times slower than the hybrid method.''}
{Yes}
{Jacobi method}
{{\cmark} Compares runtime at equal accuracy.}
{{\xmark} For Poisson's equation, Multigrid methods (or preconditioners) are highly efficient. See, e.g., figure 4e of article {\tang} or figure 23 of article \cheng. Note also that for sufficiently small 2D problems, direct solves (such as LU decomposition) will outperform Multigrid methods by a large factor.}
{N/A}
{\xmark}

\appendixentry{\stachenfeld}{\citestachenfeld}{Learned Coarse Models for Efficient Turbulence Simulation}{\stachenfeldcitations}
\appendixdata{2D incompressible \& 3D compressible Navier-Stokes}
{``Broadly, we conclude that our learned simulator outperforms traditional solvers run on coarser grids.'' Supported by figure 2 and ``Running Time'' section.}
{Yes}
{Athena++ with HLLC flux}
{{\cmark} It is fair to say that the ML-based solver outperforms the Athena++ baseline (compare the accuracy in figure 2 to the numbers in the ``running time'' section), because this comparison can be done at equal accuracy. It is not fair to say that the ML-based solver is $1000\times$ faster than Athena++ (see page 8).}
{{\dangersign} Athena++ is a state-of-the-art FV software package for shock-dominated problems.
DG and spectral methods are state-of-the-art for incompressible turbulence without shocks. In the videos available at \url{https://sites.google.com/view/learned-turbulence-simulators}, the compressible decaying turbulence seems to be weakly compressible (no shocks). See, e.g., \cite{markert2022discontinuous}, which writes that ``At least for subsonic turbulence, high order DG offers significant benefits [over FV methods] in computational efficiency for reaching a desired target accuracy.'' We give a warning sign because we believe that DG and/or spectral methods are state-of-the-art for weakly compressible turbulence and would likely outperform Athena++ with HLLC flux on this test problem, but we haven't replicated the result and so we don't have enough evidence to say for sure. }
{N/A}
{{\dangersign}}

\appendixentry{\han}{\citehan}{Predicting Physics in Mesh-reduced Space with Temporal Attention}{\hancitations}
\appendixdata{2D incompressible \& compressible Navier-Stokes}
{``We compare the evaluation cost of the learned model with the FV-based numerical models in table 6, and observe significant speedups for all three datasets.'' 100, 682, and 800 speedup reported in table 6.}
{No}
{OpenFOAM, open-source FV library.}
{{\xmark} Doesn't reduce resolution to match accuracy.}
{{\dangersign} While we recommend comparing with both FV and DG methods for the Navier-Stokes equations, we consider FV highly efficient for problems with fluid-structure interaction. However, we give a warning sign for rule 2 because OpenFOAM is known to be an inefficient general-purpose tool, with high overhead and slow convergence. See, e.g., \cite{capuanocost}. Because OpenFOAM is a general-purpose tool, we also recommend being cautious when evaluating this comparison (see recommendation 1 in Methods).}
{N/A}
{\xmark}

\appendixentry{\stevensb}{\citestevensb}
{FiniteNet: A Fully Convolutional LSTM Network Architecture for Time-Dependent Partial Differential Equations}{\stevensbcitations}
\appendixdata{1D advection, Burgers', Kuramoto-Sivashinsky}
{``We train the network on simulation data, and show that our network can reduce error by a factor of 2 to 3 compared to the baseline algorithms.'' Supported by table 1.}
{No}
{WENO5 for Burgers', 4th order finite difference for KS}
{{\cmark} Compares error (accuracy) at equal resolution (a proxy for runtime).  }
{{\cmark} WENO is highly efficient for Burgers', high-order finite difference is highly efficient for KS.}
{N/A}
{\cmark}

\appendixentry{\ozbay}{\citeozbay}
{Poisson CNN: Convolutional neural networks for the solution of the Poisson equation on a Cartesian mesh}{\ozbaycitations}
\appendixdata{2D Poisson}
{``Analytical test cases indicate that our CNN architecture is capable of predicting the correct solution of a Poisson problem with mean percentage errors below 10\%, an improvement by comparison to the first step of conventional iterative methods.'' Supported by table 6, conclusion.}
{No}
{Algebraic multigrid package PyAMG with a tolerance of $10^{-10}$}
{{\cmark} Similar speed (and theoretically faster speed on large grid sizes) and superior accuracy compared to a single cycle of multigrid. Comparing accuracy at equal runtime.}
{{\cmark} This method is focused on solving 2D Poisson problems at large grid sizes, for which multigrid methods are state-of-the-art. Note that for sufficiently small 2D problems, direct solves (such as LU decomposition) will outperform Multigrid methods. }
{N/A}
{\cmark }

\appendixentry{\zili}{\citezili}
{Graph neural network-accelerated Lagrangian fluid simulation}{\zilicitations}
\appendixdata{}
{``Overall, FGN achieves $\sim$5-8$\times$ acceleration over MPS under different resolution.'' Supported by figure 10.}
{No}
{Moving particle semi-implicit method (MPS), with conjugate gradient for pressure solve. ``MPS is a numerical method based on SPH [smooth particle hydrodynamics] which prioritizes accuracy over calculation speed.'' }
{{\cmark} ``During the benchmark, we set the absolute tolerance of CG solver to be 0.1 and maximum iteration to be 10 (note that these hyperparameters are $10^{-5}$ and 100 respectively when used to generate training data).'' When comparing runtime, they reduce the tolerance (accuracy) of the CG solver. The tolerance is chosen to be the minimum iteration that the simulation can still run without throwing NaNs (private communication with authors).
}
{{\xmark} The iterative pressure projection (i.e. Poisson solve) accounts for most of the calculation time in the baseline (see table 5). The Poisson solve uses conjugate gradient (CG), which is much slower than a state-of-the-art method such as algebraic multigrid. See, for example, \cite{sodersten2019bucket}. }
{N/A}
{\xmark}

\appendixentry{\peng}{\citepeng}{Attention-Enhanced Neural Network Models for Turbulence Simulation}{\pengcitations}
\appendixdata{2D incompressible Navier-Stokes}
{``Both models provide 8000 folds speedup compared with the pseudo-spectral numerical solver.'' Supported by table II.}
{No}
{}
{{\xmark} The ML-based solvers have non-zero error (see table 1) while the spectral solver has zero error. They are not comparing runtime at equal accuracy.}
{{\xmark} A pseudo-spectral solver for the 2D compressible Navier-Stokes with high Reynolds number on a $64\times64$ grid for 10 timesteps should take much fewer than 502 seconds to run. This is a slow implementation compared to the JAX-CFD \cite{dresdner2022learning} implementation, see \url{https://github.com/nickmcgreivy/WeakBaselinesMLPDE/blob/main/article4/data/runtime\_corr.png} where a pseudo-spectral solver on a $64\times64$ grid takes 0.1s to advance forward one unit of time.}
{N/A}
{\xmark}

\appendixentry{\chen}{\citechen}{Numerical investigation of minimum drag profiles in laminar flow using deep learning surrogates}{\chencitations}
\appendixdata{2D incompressible Navier-Stokes airfoil, steady laminar flow}
{``Therefore, relative to OpenFOAM, the speed-up factor is between 600X and 300X.'' Supported by table 2.}
{No}
{SimpleFOAM (Semi-Implicit Method for Pressure-Linked Equations, or SIMPLE), a second-order FV method within OpenFOAM}
{{\xmark} The ``performance'' section compares the runtime (for 200 iterations) of the ``ground truth'' OpenFOAM simulation to that of the ML-based surrogate solvers which are ``close to the OpenFOAM result'' but still have some error in the flow (see figure 10). To satisfy rule 1, increase the minimum mesh size of the OpenFOAM simulation to match the error of the ML-based solver, then compare runtimes at equal error.}
{{\dangersign} We give a warning sign for rule 2 because OpenFOAM is known to be an inefficient general-purpose tool, with high overhead and slow convergence. See, e.g., \cite{capuanocost}. Because OpenFOAM is a general-purpose tool, we recommend being cautious when evaluating this comparison (see recommendation 1 in Methods). }
{N/A}
{\xmark}

\appendixentry{\alguacilb}{\citealguacilb}{Predicting the propagation of acoustic waves using deep convolutional neural networks}{\alguacilbcitations}
\appendixdata{2D acoustic wave propagation, derived from Boltzmann equation (see Appendix B). Low LBM viscosity, closed domain with hard reflecting walls. ``Acoustic propagation takes place in a linear regime.'' ``The fluctuating density $\rho'$ is chosen such that $\rho_0 \gg \rho'$ to avoid non-linear effects.''}
{``The combination of both strategies can achieve a speed-up of 15.5 times with respect to the LBM code.''}
{No}
{Multi-physics lattice Boltzmann solver Palabos.}
{{\xmark} Compares the runtime of the ground truth Palabos solution to that of the ML-based solver with non-zero error.}
{{\dangersign} A Lattice Boltzmann solver is an inefficient numerical method to use to solve the linear acoustic wave equation. We give a warning sign because we believe that it would likely be much more efficient to use a linear acoustic solver or an Euler solver, but we haven't replicated this PDE and so we don't have enough evidence to say for sure. Though we don't have reason to believe that Palabos is an inefficient general-purpose tool, we also recommend being cautious when evaluating this comparison (see recommendation 1 in Methods).}
{N/A}
{{\xmark}}

\appendixentry{\wandelb}{\citewandelb}{Teaching the incompressible Navier–Stokes equations to fast neural surrogate models in three dimensions}{\wandelbcitations}
\appendixdata{3D compressible navier-stokes wake flow}
{``Furthermore, the U-Net as well as the pruned U-Net are considerably faster than Phiflow since they only require one forward pass through a convolutional neural network which can be easily parallelized and Phiflow relies on an iterative conjugate gradient solver.'' Supported by table 1.}
{No}
{Phiflow (FV method, relies on iterative conjugate gradient solver) using $128\times 128 \times 64$ grid.}
{{\xmark} The runtime of the two methods are compared in table 1, and the value of the loss function of the two methods are compared in table 1, but the accuracy of the two methods are never compared. The loss (which in this case, doesn't measure accuracy because the pressure is evolved independently of the velocity when pressure is not independent of velocity) is compared, but even so, the runtime is not compared at equal values of the loss.}
{{\cmark} While we recommend comparing with both FV and DG methods for the Navier-Stokes equations, we consider FV highly efficient for problems with fluid-structure interaction. Though we don't have reason to believe that PhiFlow is an inefficient general-purpose tool, we also recommend being cautious when evaluating this comparison (see recommendation 1 in Methods).}
{N/A}
{{\xmark}}

\appendixentry{\blist}{\citeblist}{Learned Turbulence Modelling with Differentiable Fluid Solvers: Physics-based Loss-functions and Optimisation Horizons}{\blistcitations}
\appendixdata{2D incompressible Navier-Stokes cylindrical flow}
{``For the former, our model evenly matches the performance of a $4\times$ simulation for several hundred time steps, which represents a speedup of 14.4.'' ``Measures speedups of up to 14 with respect to comparably accurate solutions from traditional solvers.'' Supported by figure 21 and table 8.}
{No}
{Semi-implicit PISO scheme, a second-order FV method}
{{\cmark} Figure 21 compares the accuracy of the ML-based solver to methods with lower resolution than ground truth. The article then compares runtime at equal accuracy.}
{{\xmark} Pseudo-spectral and DG methods are state-of-the-art for 2D incompressible Navier-Stokes. Finite-volume methods are significantly less efficient.}
{N/A}
{\xmark}

\appendixentry{\cheng}{\citecheng}{Using neural networks to solve the 2D Poisson equation for electric field computation in plasma fluid simulations}{\chengcitations}
\appendixdata{2D Poisson}
{``For this configuration, the resolution time [sic] of the neural network running on A100 GPU is about a factor 2 [sic] lower than the linear system solver on 128 cores, making it a viable option in terms of performance.'' Supported by figure 17. }
{No}
{Multigrid preconditioner with PETSc (see figure 23).}
{{\cmark} Reduces tolerance (see figure 16) to match accuracy.}
{{\cmark} Uses multigrid method. Note that for small domains, direct solves (such as LU decomposition) outperform multigrid.}
{N/A}
{\cmark}

\appendixentry{\wen}{\citewen}{An edge detector based on artificial neural network with application to hybrid compact-WENO finite difference scheme}{\wencitations}
\appendixdata{1D \& 2D shallow water and Euler}
{There are lots of results in this article, so it is hard to tell which should be considered the primary outcome. We will choose the statement in the abstract ``the ANN edge detector can capture an edge accurately with fewer grid points than the classical multi-resolution analysis.'' We also considered choosing the statement ``Generally speaking, the hybrid-ANN scheme captures the shock waves and high gradients more accurate [sic] than the Hybrid-MR scheme.'' }
{No}
{Hybrid-MR (multi-resolution analysis) scheme}
{{\cmark} The accuracy is higher (see figure 12, 13, etc) while the runtime is about the same (see table 6, 7, etc) compared to the Hybrid-MR scheme. }
{{\cmark} This is probably a strong baseline, but we are unsure.}
{}
{{\cmark}}

\appendixentry{\delaraa}{\citedelaraa}{Accelerating high order discontinuous Galerkin solvers using neural networks: 1D Burgers' equation}{\delaraacitations}
\appendixdata{1D Burgers'}
{``We see a substantial increase in efficiency, leading to ratios $CHO/CLO = 75$, 22, and 59 is case 1, 2, and 3 respectively.'' Also, ``the method shows potential and significant cost savings for high-order polynomials.'' Supported by table 1.}
{No}
{DG, very high order (polynomial 5, 7, and 28)}
{{\xmark} The ratio of cost is between a highly accurate high-order simulation and a less accurate ML-corrected low-order simulation. Thus, the cost ratio is not computed at constant accuracy.}
{{\cmark} For Burgers' equation with smooth solutions, we consider very high order DG methods highly efficient.}
{N/A}
{{\xmark} We replicate the test case (see first paragraph of section 4) using our own DG code on CPU, written in Python with JAX, and find that when the error is about $6 \times 10^{-3}$ our DG baselines take about 0.05s to run. Our baseline is 4 to 10$\times$ faster than the low-order ML-based solvers listed in table 1.}

\appendixentry{\zhao}{\citezhao}{Learning to Solve PDE-constrained Inverse Problems with Graph Networks}{\zhaocitations}
\appendixdata{2D scalar wave equation, 2D incompressible Navier-Stokes}
{``The proposed method \dots is $35\times$ faster than the classical FEM solver.'' Supported by table 2.}
{No}
{Wave equation: FEniCS using Euler Method and first order elements, with GMRES and LU factorization as preconditioner. Navier-stokes: Chorin's method.}
{{\dangersign} Table 2 compares the runtime at equal resolution (a proxy for runtime) rather than at equal accuracy. We give this article a warning sign because while it does reduce the resolution (relative to the fine grid ground truth FEM solver) it doesn't reduce it enough to reach equal accuracy. The ML-based solver has an error $2.5\times$ higher than the FEM solver at equal resolution, and thus we think this comparison is likely unfair.}
{{\cmark} While we recommend comparing with both FV and DG methods for the Navier-Stokes equations, we consider FV highly efficient for problems with fluid-structure interaction. Though we don't have reason to believe that FEniCS is an inefficient general-purpose tool, we also recommend being cautious when evaluating this comparison (see recommendation 1 in Methods).}
{N/A}
{{\dangersign}}

\appendixentry{\assessments}{\citeassessments}{Performance and accuracy assessments of an incompressible fluid solver coupled with a deep Convolutional Neural Network}{\assessmentscitations}
\appendixdata{2D Poisson}
{``These networks can provide solutions 10-25 faster than traditional iterative solvers.''}
{No}
{Jacobi method}
{{\cmark} ``A fair performance comparisons is performed \dots allowing the assessment of the time of inference at a fixed error level.'' }
{{\xmark} A multigrid method or preconditioner is considered highly efficient for elliptic PDEs.}
{N/A}
{\xmark}

\appendixentry{\holloway}{\citeholloway}{Acceleration of Boltzmann Collision Integral Calculation Using Machine Learning}{\hollowaycitations}
\appendixdata{6-dimensional Boltzmann collision operator}
{``Our method demonstrated a speed up of 270 times compared to these methods while still maintaining reasonable accuracy.'' Supported by table 1.}
{No}
{DG discretization of collision operator}
{{\xmark} Compares runtime of highly accurate ground DG method to that of less accurate ML-based solver.}
{{\cmark} We believe that the DG discretization is likely highly efficient.}
{N/A}
{\xmark}

\appendixentry{\azulay}{\citeazulay}{Multigrid-augmented deep learning preconditioners for the Helmholtz equation}{\azulaycitations}
\appendixdata{2D Helmholtz}
{``We show that while our U-Net may require more FLOPs than traditional methods, it can applied efficiently on GPU hardware, and yield favorable running times.'' Supported by figure 10 and section 4.3.5.}
{No}
{Geometric multigrid preconditioner, followed by GMRES iterations }
{{\cmark} Figure 10 compares runtime at constant accuracy/error.}
{{\cmark} Multigrid preconditioner is highly efficient for elliptic PDEs.}
{}
{\cmark}

\appendixentry{\wulatent}{\citewulatent}{Learning to Accelerate Partial Differential Equations via Latent Global Evolution}{\wulatentcitations}
\appendixdata{1D Burgers', 2D incompressible Navier-Stokes, 3D incompressible Navier-Stokes cylinder flow }
{``It achieves significant [sic] smaller runtime compared to the MP-PDE model (which is much faster than the classical WENO5 scheme).'' Also, ``we see that our LE-PDE achieves a $70.80/0.084 \simeq 840\times$ speed up compared to the ground-truth solver on the same GPU.'' See supplementary material, figure S1, and table 5.}
{No}
{PhiFlow}
{{\xmark} Table 5 compares the runtime of a highly accurate PhiFlow to that of a less accurate ML-based solver.}
{{\cmark} While we recommend comparing with both FV and DG methods for the Navier-Stokes equations, we consider FV highly efficient for problems with fluid-structure interaction. Though we don't have reason to believe that PhiFlow is an inefficient general-purpose tool, we also recommend being cautious when evaluating this comparison (see recommendation 1 in Methods).}
{N/A}
{\xmark}

\appendixentry{\liu}{\citeliu}{Predicting parametric spatiotemporal dynamics by multi-resolution PDE structure-preserved deep learning}{\liucitations}
\appendixdata{2D Burgers', 2D incompressible Navier-Stokes}
{``In particular, the speedup by the PPNN varies from $10\times$ to $60\times$ without notably sacrificing the prediction accuracy.'' Supported by figure 7b.}
{No}
{Burgers': 3rd-order accurate up-wind scheme for convection, 6th order accurate central-difference scheme for diffusion term, forward-euler timestepping. Navier-stokes: PISO algorithm using OpenFOAM.}
{{\xmark} Doesn't reduce the resolution of the PISO algorithm to compare at equal accuracy.}
{{\xmark} Pseudo-spectral or DG algorithms are considered highly efficient for 2D incompressible Navier-stokes. Finite-volume algorithms like PISO are much less efficient. Furthermore, OpenFOAM is known to be an inefficient general-purpose tool, with high overhead and slow convergence. See, e.g., \cite{capuanocost}. Because OpenFOAM is a general purpose tool, we also recommend being cautious when evaluating this comparison (see recommendation 1 in Methods). }
{N/A}
{\xmark}

\appendixentry{\zhang}{\citezhang}{A Hybrid Iterative Numerical Transferable Solver (HINTS) for PDEs Based on Deep Operator Network and Relaxation Methods}{\zhangcitations}
\appendixdata{1D \& 2D Poisson, 1D 2D \& 3D Helmholtz}
{``The results show that HINTS performs consistently better than the corresponding numerical methods, with the improvement of computational efficiency being up to $\mathcal{O}(10^2)$.'' Supported by figure S9.}
{No}
{For 1D Poisson: Multigrid with damped Jacobi relaxation, either 3 or 5 grid levels.}
{{\cmark} Figure S9 measures time until convergence, which implies equal accuracy.}
{{\xmark} For small problems (such as 1D Poisson), direct methods like LU decomposition are much faster than iterative methods like multigrid. A 1D Poisson problem should take microseconds to milliseconds to solve to machine precision, not 1 to 100 seconds (as in figure S9).}
{N/A}
{\xmark}

\appendixentry{\duarte}{\citeduarte}{Black hole weather forecasting with deep learning: a pilot study}{\duartecitations}
\appendixdata{2D black hole hydrodynamics: ``Our data set was generated from two-dimensional hydrodynamical simulations of viscous accretion on to a Schwarzschild BH.''}
{``For instance, once trained the model evolves an RIAF on a single GPU four orders of magnitude faster than usual fluid dynamics integrators running in parallel on 200 CPU cores.'' Supported by table 3, but caveated by the comments in section 5.1.}
{Yes}
{PLUTO code which uses Godunov-like flux}
{{\xmark} This article compares the runtime of a highly accurate standard numerical method to that of a less accurate ML-based solver. Doesn't compare runtime at equal accuracy.}
{{\cmark} We consider FV schemes highly efficient for this problem.}
{}
{\xmark}

\appendixentry{\alguacil}{\citealguacil}{Deep Learning Surrogate for the Temporal Propagation and Scattering of Acoustic Waves}{\alguacilcitations}
\appendixdata{2D (acoustic) wave}
{``When both strategies are combined, a large acceleration factor of 141 can be achieved with respect to the MPI-based simulation.'' Supported by table 5.}
{No}
{Lattice-Boltzmann simulation Palabos, same as article \alguacilb.}
{{\xmark} Table 5 compares the runtime of the highly accurate baseline to that of the less accurate ML-based solver.}
{{\dangersign} A Lattice Boltzmann solver is an inefficient numerical method to use to solve the linear acoustic wave equation. We give a warning sign because we believe that it would likely be much more efficient to use a linear acoustic solver or an Euler solver, but we haven't replicated this PDE and so we don't have enough evidence to say for sure. Though we don't have reason to believe that Palabos is an inefficient general-purpose tool, we also recommend being cautious when evaluating this comparison (see recommendation 1 in Methods).}
{N/A}
{\xmark}

\appendixentry{\bezginb}{\citebezginb}{WENO3-NN: A maximum-order three-point data-driven weighted essentially non-oscillatory scheme}{\bezginbcitations}
\appendixdata{1D advection, 1D \& 2D Euler}
{``The WENO3-NN scheme shows very good generalizability across all benchmark cases and different resolutions, and exhibits a performance similar to or better than the classical WENO5-JS scheme.'' Supported by table 4.}
{No}
{WENO5-JS}
{{\cmark} Table 4 and various figures compares accuracy at constant resolution, a proxy for runtime.}
{{\xmark} The primary outcome (in this case defined by what is reported in the abstract) compares the ML-based solver (WENO3-NN) to a baseline (WENO5-JS). As we learn in section 4.4 and in the appendix, both WENO5-JS and WENO3-NN underperform relative to WENO5-Z for 5 or 6 out of the 6 benchmark problems. Because the primary outcome compares to the weaker baseline (WENO5-JS) rather than the stronger baseline (WENO5-Z), we determine that rule 2 is not satisfied.
}
{N/A}
{\xmark }

\appendixentry{\shang}{\citeshang}{Deep Petrov-Galerkin Method for Solving Partial Differential Equations}{\shangcitations}
\appendixdata{2D Poisson, 2D Wave}
{``This new method outperforms traditional numerical methods in several aspects: compared to the finite element method and finite difference method, DPGM is much more accurate with respect to degrees of freedom.'' Supported by tables 6 and 7.}
{No}
{FEniCS, lagrange elements $P_k$ $k=1,2,3$ for Poisson and }
{{\cmark} Compares accuracy at constant degrees of freedom (which is a proxy for runtime). }
{{\cmark} Though we don't have reason to believe that FEniCS is an inefficient general-purpose tool for these PDEs, we recommend being cautious when evaluating this comparison (see recommendation 1 in Methods).}
{N/A}
{\cmark}

\appendixentry{\kube}{\citekube}{Machine learning accelerated particle-in-cell plasma simulations}{\kubecitations}
\appendixdata{1D particle-in-cell (linear solver)}
{``We find that this approach reduces the average number of required solver iterations by about 25\% when simulating electron plasma oscillations.'' Supported by figure 2.}
{No}
{GMRES (Jacobian-free Newton-Krylov method) to solve system of non-linear equations}
{{\xmark} Although this article compares the number of iterations at equal accuracy, (based on private communication with the first author) the runtime of the ML-based solver is significantly longer than the standard GMRES solver. We gave an `X' because, in this case, the number of iterations is not a fair proxy for speed.}
{{\cmark} Unsure. }
{N/A}
{\xmark}

\appendixentry{\shi}{\citeshi}{LordNet: Learning to Solve Parametric Partial Differential Equations without Simulated Data}{\shicitations}
\appendixdata{2D Poisson, 2D incompressible Navier-Stokes}
{``For Navier-Stokes equation, the learned operator is over 50 times faster than the finite difference solution with the same computational resources.'' Supported by figure 3 and section 4.2.}
{No}
{Finite-difference scheme (FDM) with central differencing, with conjugate gradient method to solve sparse algebraic equations}
{{\xmark} Doesn't reduce resolution of FDM to match accuracy of ML-based solver.}
{{\cmark} Unsure.}
{}
{\xmark}

\appendixentry{\ranadea}{\citeranadea}{A Latent space solver for PDE generalization}{\ranadeacitations}
\appendixdata{3D steady-state electronic cooling with natural convection}
{``Thus, the hybrid solver results in a 200$\times$ speedup over Ansys Fluent in generating solutions on fine girds [sic].'' Supported by section 3.2.1.}
{No}
{ANSYS Fluent}
{{\xmark} Doesn't reduce resolution of ANSYS baseline to match accuracy of ML-based solver. }
{{\cmark} Though we don't have reason to believe that ANSYS is an inefficient general-purpose tool, we recommend being cautious when evaluating this comparison (see recommendation 1 in Methods).}
{N/A}
{\xmark}

\appendixentry{\chenb}{\citechenb}{A machine learning based solver for pressure Poisson equations}{\chenbcitations}
\appendixdata{2D Poisson}
{``The ML-block provides a better initial iteration value for the traditional iterative solver, which greatly reduces the number of iterations of the traditional iterative solver and speeds up the solution of the PPE.'' Supported by table 2 and equation 15.}
{No}
{Preconditioned conjugate gradient}
{{\cmark} Equation 15: ``to achieve the same solution accuracy''. Compares number of iterations at equal accuracy. }
{{\xmark} Multigrid methods (or preconditioners) are state-of-the-art for Poisson's equation. See, e.g., figure 4e of article {\tang} or figure 23 of article \cheng. Note also that for sufficiently small 2D problems, direct solves (such as LU decomposition) will outperform Multigrid methods.}
{N/A}
{\xmark}

\appendixentry{\ranadeb}{\citeranadeb}{A composable autoencoder-based iterative algorithm for accelerating numerical simulations}{\ranadebcitations}
\appendixdata{2D Laplace, 2D incompressible Navier-Stokes, 3D electronic cooling with natural convection, 3D steady-state channel flow (incompressible Navier-Stokes?), 3D transient flow over a cylinder with heating}
{``We observe that the CoAE-MLSim approach is about 40-50$\times$ faster in the steady-state cases and about 100$\times$ faster in transient cases as compared to commercial PDE solvers such as Ansys Fluent for the experiments presented in this work.'' Not supported by any quantitative analysis.}
{No}
{Ansys Fluent}
{{\xmark} Doesn't compare the runtime at equal accuracy, instead compares runtime of highly accurate ANSYS Fluent to less accurate ML-based solver.}
{{\cmark} Though we don't have reason to believe that Ansys is an inefficient general-purpose tool, we recommend being cautious when evaluating this comparison (see recommendation 1 in Methods).}
{N/A}
{{\xmark}}

\appendixentry{\pengb}{\citepengb}{Linear attention coupled Fourier neural operator for simulation of three-dimensional turbulence}{\pengbcitations}
\appendixdata{3D incompressible Navier-Stokes}
{``During inference, both neural network models provides [sic] 20 folds speedup compared with the DNS approach with the traditional numerical solver.'' Supported by table 4.}
{No}
{A pseudo-spectral method, $64\times 64\times 64$ box.}
{{\xmark} This article never reduces the resolution of the spectral baseline to compare at equal accuracy.}
{{\cmark} Pseudo-spectral methods are highly efficient for this problem.}
{N/A}
{\xmark}

\appendixentry{\delarab}{\citedelarab}{Accelerating high order discontinuous Galerkin solvers using neural networks: 3D compressible Navier-Stokes equations}{\delarabcitations}
\appendixdata{3D compressible Navier-Stokes in the incompressible limit}
{``The low order corrected solution is 4 to 5 times faster than a simulation with comparable accuracy (polynomial order 6).'' Supported by section 3.5, figure 12, table 2.}
{Yes}
{Very high order DG (polyOrder 5 to 6) with constant number of elements. HORSES3D, a nodal high-order DG spectral element method (DGSEM).}
{{\cmark} Compares runtime at equal accuracy (see section 3.5).}
{{\cmark} We consider DG methods highly efficient for weakly compressible flows.
}
{N/A}
{\cmark}

\appendixentry{\ranadec}{\citeranadec}{A composable machine-learning approach for steady-state simulations on high-resolution grids}{\ranadeccitations}
\appendixdata{2D Laplace, 2D Poisson, 2D non-linear Poisson, 3D electronic cooling with natural convection}
{``We observe that the CoAE-MLSim approach is about 40-50$\times$ faster as compared to commercial steady-state PDE solvers such as ANSYS Fluent for the same mesh resolution and physical domain size in all the experiments presented in this work.'' See appendix E. }
{Yes}
{Ansys Fluent}
{{\xmark} Compares runtime at constant mesh resolution (a proxy for runtime), not at constant accuracy. For a fair comparison, reduce the resolution of ANSYS Fluent to match the accuracy of the ML-based solver.}
{{\cmark} Though we don't have reason to believe that Ansys is an inefficient general-purpose tool, we recommend being cautious when evaluating this comparison (see recommendation 1 in Methods).}
{N/A}
{\xmark}

\appendixentry{\fang}{\citefang}{Immersed boundary-physics informed machine learning approach for fluid–solid coupling}{\fangcitations}
\appendixdata{2D incompressible Navier-Stokes cylindrical flow, oscillating cylinder}
{``The time consumed by the machine learning model is reduced by 38.5\% compared with IB-LBM.'' Supported by figure 12.}
{No}
{Immersed-boundary lattice boltzmann method (IB-LBM). A D2Q9 scheme (see section 2.1).  }
{{\xmark} This comparison is not done at equal accuracy. See the differences between figure 5 and figure 9, as well as table 4. Since this paper is using the drag coefficient as a measure of accuracy, a fair comparison of simulation time would be at equal drag coefficient (or equal deviation from the high-resolution simulation's drag coefficient). }
{{\cmark} Unsure.}
{N/A}
{\xmark}

\appendixentry{\shukla}{\citeshukla}{Deep neural operators can serve as accurate surrogates for shape optimization: A case study for airfoils}{\shuklacitations}
\appendixdata{2D compressible Navier-Stokes airfoil wing}
{Importantly, DeepONets exhibit almost no generalization error over the dataset, so it follows that the resulting optimized geometry is accurate and achieved 32, 253 speed-up compared to the CFD baseline.'' Supported by table 4.}
{No}
{Nektar, using DG spectral element method (DGSEM) with basis functions spanned in 2D by Legendre polynomials of the second degree. DIRK used for time integration. See section 3.1.3.}
{{\xmark} While the Nektar solution has zero error, the DeepONet surrogate model has positive error: ``a closer look reveals non-physical streamlines originating from the surface of the airfoil. This is due to the error in the velocity fields near the airfoil surface, predicted by the DeepONet.'' Thus, the comparison in table 4 is not done at constant accuracy. }
{{\cmark} DGSEM is likely highly efficient for this problem. Though we don't have reason to believe that Nektar is an inefficient general-purpose tool, we recommend being cautious when evaluating this comparison (see recommendation 1 in Methods).}
{N/A}
{\xmark}

\appendixentry{\zhangb}{\citezhangb}{Learning the elastic wave equation with Fourier Neural Operators}{\zhangbcitations}
\appendixdata{(elastic) wave equation}
{``Post-training, the FNO is observed to generate accurate elastic wave fields at approximately 10 times the speed of traditional finite difference methods.'' Supported by table 1 and conclusion.}
{No}
{``Wavefields generated with the isotropic stress-velocity wave equation, using a staggered grid finite difference method with 4th order accuracy in space as training data.'' 84 by 84 grid. }
{{\xmark} Compares the runtime of the ground truth finite difference method to that of the less accurate ML-based solver. Doesn't reduce grid resolution to compare at equal accuracy. }
{{\cmark} Unsure.}
{N/A}
{\xmark }

\appendixentry{\bezgin}{\citebezgin}{A fully-differentiable compressible high-order computational fluid dynamics solver}{\bezgincitations}
\appendixdata{2D \& 3D compressible Navier-Stokes}
{``The NN-Rusanov flux stays stable over the course of the simulation and consistently outperforms the Rusanov flux\dots The NN-Rusanov flux is less dissipative than the classical Rusanov scheme.'' Supported by figure 5.}
{No}
{Rusanov flux}
{{\cmark} Compares error (accuracy) at equal resolution (a proxy for runtime) in figure 5. }
{{\xmark} Rusanov flux is a first-order method (with ``excess numerical diffusion [that] leads to a smeared out solution'') and is not a highly efficient baseline for this problem. }
{N/A}
{\xmark}

\appendixentry{\yang}{\citeyang}{Rapid Seismic Waveform Modeling and Inversion with Neural Operators}{\yangcitations}
\appendixdata{2D acoustic wave equation}
{``We show that full waveform modeling with neural operators is nearly two orders of magnitude faster than conventional numerical methods.'' Supported by page 2: runtime of 0.02 sec versus 1.23 sec. }
{No}
{Finite difference code on a $64 \times 64$ mesh.}
{{\xmark} Doesn't reduce resolution to compare at equal accuracy. }
{{\cmark} The finite-difference baseline in reference 40 seems likely to be a strong baseline, but we aren't entirely sure.}
{N/A}
{\xmark}

\appendixentry{\tang}{\citetang}{Neural Green’s function for Laplacian systems}{\tangcitations}
\appendixdata{2D Poisson}
{``Although our method saves only about $2\times$ the number of multiply-add operations compared to MGPCG, its intrinsic parallel nature enables it to reach a speedup of up to $12\times$ at all resolutions.'' Supported by table 2 and figure 4b. }
{No}
{Multigrid and multigrid-preconditioned conjugate gradient (MGPCG) are the two strongest baselines. }
{{\cmark} Compares at equal accuracy in table 2.}
{{\xmark} Consider table 2: the strongest baseline takes 90ms, 189ms, 299ms, and 425ms for grid sizes of $33\times33$, $65\times65$, $129\times129$, and $257\times257$ respectively. We implement Poisson's equation on a square domain using a direct solve (LU decomposition) and find that the direct solve takes 0.2ms, 0.4ms, 3ms, and 12ms for grid sizes of $32\times32$, $64\times64$, $128\times128$, and $256\times256$ respectively. LU decomposition is between 500 and 35 times faster than multigrid. Multigrid a weak baseline relative to direct methods for sufficiently small problems.}
{N/A}
{\xmark}

\appendixentry{\nastorg}{\citenastorg}{DS-GPS : A Deep Statistical Graph Poisson Solver (for faster CFD simulations)}{\nastorgcitations}
\appendixdata{2D poisson}
{``By taking advantage of GPU parallelism, we observe that our method can compute the solution ten times faster than LU decomposition.'' Not supported by any evidence.}
{No}
{LU decomposition}
{{\xmark} The error of LU decomposition on the graph is zero. The error of this ML-based solver is higher than zero (see figure 1). For a fair comparison, compare the runtime of LU decomposition on a coarsened graph (so that the error on the fine graph is comparable).}
{{\cmark} For small problems, LU decomposition is highly efficient or state-of-the-art.}
{N/A}
{\xmark}

\appendixentry{\gopakumar}{\citegopakumar}{Fourier Neural Operator for Plasma Modelling}{\gopakumarcitations}
\appendixdata{2D magnetohydrodynamics (MHD)}
{``Our work shows that the FNO is capable of predicting magnetohy- drodynamic models governing the plasma dynamics, 6 orders of magnitude faster than the traditional numerical solver, while maintaining considerable accuracy (NMSE $\sim10^5$).'' Supported by table 1.}
{No}
{JOREK code. 200 by 200 bi-cubic finite-elements.}
{{\xmark} Doesn't reduce accuracy of JOREK to match accuracy of FNO. See errors in, e.g., figure 2.}
{{\dangersign} 160 CPU hours to solve a 2D advection-dominated time-dependent PDE is unusually long, even at a high resolution of $200\times200$. We believe that JOREK is likely using inefficient numerical methods for advection-dominated flows. We believe the most likely explanation for the slow implementation is that JOREK uses fully implicit timestepping, while advection-dominated PDEs are better suited for explicit timestepping. 
}
{N/A}
{\xmark}

\appendixentry{\shit}{\citeshit}{Semi-Implicit Neural Solver for Time-dependent Partial Differential Equations}{\shitcitations}
\appendixdata{2D advection-diffusion, with advection speeds in range [-2, 2] and diffusion coefficients in range [0.2, 0.8]. }
{``We observe that for a given acceptance, the neural solver is 19.2\% faster.'' Supported by figure 5.}
{No}
{Semi-implicit scheme (see equation 3), using fixed-point iteration introduced in section 2.1.}
{{\cmark} Figure 5 compares runtime at equal error.}
{{\xmark} For this mixed hyperbolic-parabolic PDE, an (explicit) super-time-stepping method (e.g., \cite{meyer2014stabilized}) would be more efficient than an implicit method. Let us explain why. The semi-implicit scheme uses a timestep of $\Delta t_{\textnormal{implicit}} = 0.2$, and each timestep requires a number of fixed point iterations (up to 25). An explicit method would (since the maximum diffusion coefficient is $\nu = 0.8$ and $\Delta x = 0.098$) have a timestep of $\Delta t_{\textnormal{explicit}} \approx \frac{(\Delta x)^2}{2\nu} = 0.006$. Since $\frac{\Delta t_{\textnormal{implicit}}}{\Delta t_{\textnormal{explicit}}} \approx 33$, but each implicit timestep takes roughly 25$\times$ more runtime, then implicit timestepping would likely be as fast or slightly faster (as a rough estimate, as much as $33/25\approx 1.3\times$ faster) than naive explicit timestepping. However, $s$-stage super-time-stepping allows for timesteps proportional to $s^2$ and a runtime speedup (relative to explicit timestepping) of $s^2/s = s$, where the optimal number of stages $s = \sqrt{\frac{\Delta t_{\textnormal{implicit}}}{\Delta t_{\textnormal{explicit}}}}$. In this case, the speedup from using super-time-stepping would be $s = \sqrt{33} = 5.75$, which would be faster than the implicit method. }
{N/A}
{\xmark}

\appendixentry{\suforecast}{\citesuforecast}{Forecasting Variable-Density 3D Turbulent Flow}{\suforecastcitations}
\appendixdata{3D compressible Navier-Stokes}
{``Across flows with different density-ratio, our method over [sic] 3 orders of magnitude faster than high-fidelity numerical simulations.'' Supported by table 1. }
{No}
{Petascale variable-density version of the CFDNS code. Spectral code: spatial derivatives are evaluated using Fourier transforms. Uses ($64^3$) resolution.}
{{\xmark} Table 1 compares the runtime of the ground truth solver to that of the less accurate ML-based solver. For a fair comparison, reduce the resolution of the CFDNS code until the accuracy matches the ML-based solver.}
{{\cmark} Spectral codes are highly efficient for weakly compressible turbulence.}
{N/A}
{\xmark}

\appendixentry{\jeon}{\citejeon}{Physics-Informed Transfer Learning Strategy to Accelerate Unsteady Fluid Flow Simulations}{\jeoncitations}
\appendixdata{2D incompressible Navier-Stokes}
{``The simulation was accelerated by 1.8 times in the laminar counterflow CFD dataset condition including the parameter-updating time.''}
{No}
{IcoFoam, which is part of OpenFOAM, which uses the PISO algorithm.}
{{\xmark} The abstract writes that ``Notably, the cross-coupling strategy with a grid-based network model does not compromise the simulation accuracy for computational acceleration.'' If this were true, then rule 1 would be satisfied. The problem is that it is the residual of the governing equations that is ``not compromise[d]'', not the accuracy of the solution. As the paper correctly points out at the top of page 3, a low residual does not imply low error. Yes, the cross-coupling ML-CFD strategy is $1.8\times$ faster than the ``ground truth'' CFD solver, but as figures 9 and 12 show the ``ground truth'' CFD solver has zero error while the cross-coupling ML-CFD strategy has positive error. Thus, rule 1 is not satisfied. To satisfy rule 1, increase the error of the CFD solver to match that of the cross-coupling strategy. One way to do this would be to reduce the resolution of the OpenFOAM solution; a second would be to reduce the residual of OpenFOAM below the default residual.
}
{{\dangersign} We give a warning sign for rule 2 because OpenFOAM is known to be an inefficient general-purpose tool, with high overhead and slow convergence. See, e.g., \cite{capuanocost}. Because OpenFOAM is a general purpose tool, we also recommend being cautious when evaluating this comparison (see recommendation 1 in Methods). While we recommend comparing with both FV and DG methods for the Navier-Stokes equations, we consider FV highly efficient for problems with fluid-structure interaction. }
{N/A}
{\xmark}

\appendixentry{\dai}{\citedai}{FourNetFlows: An efficient model for steady airfoil flows prediction}{\daicitations}
\appendixdata{2D incompressible Navier-Stokes airfoil wing}
{``We quantitatively shows [sic] the accuracy of FourNetFlows is matched with the traditional method, running four to five orders of magnitude faster.'' Supported by table 1. }
{No}
{SimpleFoam, which uses SIMPLE (semi-implicit method for pressure linked equations) and is part of OpenFOAM.}
{{\xmark} Compares the runtime of ground truth SimpleFOAM to that of the less accurate ML-based solver.}
{{\dangersign} We give a warning sign for rule 2 because OpenFOAM is known to be an inefficient general-purpose tool, with high overhead and slow convergence. See, e.g., \cite{capuanocost}. Because OpenFOAM is a general purpose tool, we also recommend being cautious when evaluating this comparison (see recommendation 1 in Methods).}
{N/A}
{\xmark}

\appendixentry{\sun}{\citesun}{Local Randomized Neural Networks with Discontinuous Galerkin Methods for Partial Differential Equations}{\suncitations}
\appendixdata{1D Helmholtz, 2D Poisson.}
{``We compare the proposed methods with the finite element method and the usual DG method. The LRNN-DG methods can achieve better accuracy under the same degrees of freedom.'' Supported by figure 4.}
{No}
{FEM and DG, up to 3rd order. Doesn't have any other details.}
{{\cmark} Compares accuracy at constant number of degrees of freedom, a proxy for runtime. }
{{\cmark} FEM and DG are efficient methods for solving Poisson's equation.}
{N/A}
{\cmark}

\appendixentry{\shao}{\citeshao}{A Poisson’s Equation Solver Based on Neural Network Precondtioned CG Method}{\shaocitations}
\appendixdata{2D Poisson}
{``Numerical examples demonstrate that compared to conjugate gradient (CG) method, NN-PCG significantly improves convergence performance on solving 2-D Poisson’s equation.''}
{No}
{Conjugate gradient method.}
{{\cmark} Figure 2 compares accuracy at constant number of iterations.}
{{\xmark} Multigrid methods are state-of-the-art for this PDE. See, e.g., figure 4e of article {\tang} or figure 23 of article \cheng. Note also that for sufficiently small 2D problems, direct solves (such as LU decomposition) will outperform Multigrid methods.}
{N/A}
{\xmark}

\appendixentry{\discacciati}{\citediscacciati}{Controlling oscillations in high-order Discontinuous Galerkin schemes using artificial viscosity tuned by neural networks}{\discacciaticitations}
\appendixdata{1D \& 2D advection, 1D \& 2D Burgers', 1D \& 2D Euler}
{``The network-based model is always at par with the best among the traditional optimized models.'' Supported by section 6 and tables 13 and 14.}
{No}
{N/A (doesn't claim superiority)}
{N/A}
{N/A}
{N/A}
{N/A}

\appendixentry{\magiera}{\citemagiera}{Constraint-aware neural networks for Riemann problems}{\magieracitations}
\appendixdata{1D Euler, 1D scalar hyperbolic PDE}
{``Naturally, the case study problems considered in Section 6 do not show a speedup of the computational time by exchanging the analytical Riemann solver by a surrogate neural network.''}
{No}
{N/A}
{N/A}
{N/A}
{N/A}
{N/A}

\appendixentry{\bezginc}{\citebezginc}{A data-driven physics-informed finite-volume scheme for nonclassical undercompressive shocks}{\bezginccitations}
\appendixdata{Cubic scalar hyperbolic conservation law}
{``For the weak shock test case which is calculated on a coarse mesh, the NN scheme is roughly 10 times slower than the WCD scheme.''}
{No}
{WCD scheme}
{N/A (doesn't claim superiority)}
{N/A}
{N/A}
{N/A}

\appendixentry{\dongb}{\citedongb}{On computing the hyperparameter of extreme learning machines: Algorithm and application to computational PDEs, and comparison with classical and high-order finite elements}{\dongbcitations}
\appendixdata{2D Poisson, 2D non-linear Helmholtz, 1D Burgers.}
{``It is shown that the current improved ELM far outperforms the classical FEM. Its computational performance is comparable to that of the high-order FEM for smaller problem sizes, and for larger problem sizes the ELM markedly outperforms the high-order FEM.'' Supported by figures 16, 27, and 36.}
{No}
{FEniCS}
{N/A (doesn't claim superiority to stronger baseline)}
{N/A}
{N/A}
{N/A}

\appendixentry{\dresdner}{\citedresdner \hspace{0.1cm}(Note: we consider version 1 of this article on ArXiv, version 2 was uploaded after our systematic review was completed but before this article was finished.)}{Learning to correct spectral methods for simulating turbulent flows}{\dresdnercitations}
\appendixdata{1D Kuramoto-Sivashinsky (KS), 1D Burgers', 2D incompressible Navier-Stokes }
{``Overall there is little potential for accelerating 2D turbulence beyond traditional spectral solvers.''}
{No}
{Spectral solver}
{N/A (doesn't claim superiority)}
{N/A}
{N/A}
{N/A}

\appendixentry{\toshev}{\citetoshev}{E(3) Equivariant Graph Neural Networks for Particle-Based Fluid Mechanics}{\toshevcitations}
\appendixdata{3D incompressible Navier-Stokes}
{``Our main findings are that while being rather slow to train and evaluate\dots'' Supported by table 1.}
{Yes}
{N/A}
{N/A (Doesn't claim superiority)}
{N/A}
{N/A}
{N/A}

\newpage 

%\bibliography{sn-bibliography}

\end{document}